\newif\ifshort
\shortfalse 

\documentclass{amsart} 

\ifshort
\else
\fi

\usepackage{va, stycb}
\usepackage{todonotes}

\usepackage[pdftex]{hyperref}
\hypersetup{colorlinks,citecolor=blue,linktocpage,hyperindex=true,backref=true}

\usepackage{float}
\usepackage{enumitem}
\usepackage{bm}
\usepackage{verbatim}

\let\oldtocsection=\tocsection
\let\oldtocsubsection=\tocsubsection
\let\oldtocsubsubsection=\tocsubsubsection
\renewcommand{\tocsection}[2]{\hspace{0em}\oldtocsection{#1}{#2}}
\renewcommand{\tocsubsection}[2]{\hspace{1em}\oldtocsubsection{#1}{#2}}
\renewcommand{\tocsubsubsection}[2]{\hspace{0.5em}\oldtocsubsubsection{#1}{#2}}

\title[Compact moduli spaces of Campedelli and Burniat surfaces]%
{Explicit compactifications of moduli spaces\\
  of Campedelli and Burniat surfaces}

\author{Valery Alexeev} 
\email{valery@uga.edu}
\address{Department of Mathematics, University of Georgia, Athens GA 30602, USA}
\author{Rita Pardini} 
\email{rita.pardini@unipi.it}
\address{Dipartimento di Matematica, Universit\`a di Pisa, Largo B. Pontecorvo 5, 56127 Pisa, Italy}

\begin{document}
\maketitle

\begin{abstract}
  We describe explicitly the geometric compactifications, obtained
  by adding slc surfaces~$X$ with ample canonical class, for two connected
  components in the moduli space of surfaces of general type:
  Campedelli surfaces with $\pi_1(X)=\bZ_2^3$ and Burniat
  surfaces with $K^2=6$.
\end{abstract}

\ifshort
\setcounter{tocdepth}{1}
\else
\setcounter{tocdepth}{2}
\fi

\tableofcontents

\section{Introduction}
\label{sec:intro}

In 1988, Koll\'ar and Shepherd-Barron
\cite{kollar1988threefolds-and-deformations} proposed a way to
compactify the moduli space of surfaces of general type by adding
\emph{stable surfaces}, i.e. surfaces that have slc (semi log
canonical) singularities and ample canonical class $K_X$, similar to
the stable curves in dimension one. This construction was subsequently
extended to pairs $(X,B=\sum b_iB_i)$ with divisors and to higher
dimensions.  The resulting compact moduli spaces are commonly known as
KSBA spaces, see \cite[Ch.~8]{kollar2023families-of-varieties}, and
they will be denoted $\oM\slc$ here.

In the cases with nonzero boundary divisor $B$
there are many papers
where these compact moduli spaces are described in detail, e.g. for
toric and abelian varieties \cite{alexeev2002complete-moduli},
hyperplane arrangements \cite{hacking2006compactification-moduli,
  alexeev2015moduli-weighted}, K3 surfaces \cite{alexeev2023compact,
  alexeev2023stable-pair, alexeev2022compactifications-moduli,
  alexeev2022mirror-symmetric}, elliptic surfaces
\cite{ascher2021moduli-of-weighted, inchiostro20moduli-weierstrass}
and many, many more.

However, in the original case of
\cite{kollar1988threefolds-and-deformations}, i.e. with $B=0$,
practically no explicit compactifications are known, aside from the moduli of
surfaces which are either quotients or special covers of a product of
two curves \cite{vanopstall2006stable-degenerations1,
  liu2012stable-degenerations, rollenske2010compact-moduli}, which
essentially reduce to moduli of curves.
Let us also mention some works on \emph{partial} compactifications:
\cite{Franciosi2017K21, Fantechi2022semismooth,
  Franciosi2022Tsing, coughlan2023Tdivisors, gallardo2022unimodal}.

The goal of this paper is to describe explicitly two complete
compactifications in the case with zero boundary divisor, for two
irreducible components of the classical moduli space of surfaces of
general type, of dimensions $6$ and $4$:
\begin{enumerate}
\item Campedelli surfaces with $\pi_1(X)=\bZ_2^3$. They can be defined
  as $\bZ_2^3$-covers of $\bP^2$ ramified in $7$ lines $B_i$. The
  problem then can be reduced to compactifying the moduli of pairs
  $(\bP^2,\sum_{i=1}^7 \frac12 B_i)$---which turns out to be quite easy
  in this case---and applying the theory of singular abelian covers of
  \cite{alexeev2012non-normal-abelian}.
\item Burniat surfaces with $K_X^2=6$. They can be defined as
  $\bZ_2^2$-covers of the del Pezzo surface $\Sigma=\Bl_3\bP^2$, ramified in
  $12$ curves coming from a particular configuration of $9$ lines in
  $\bP^2$. This case, although similar in spirit to the one above, turns
  out to be much harder. 
\end{enumerate}

Our main results are:

\begin{theorem-intro}\label{thm-intro:campedelli}
  For Campedelli surfaces,
  the main irreducible component of 
  the compactification $\oM\cam\slc$ is 
  \begin{displaymath}
    \GL(3,\bF_2) \backslash \big( \bP^2\big)^7//\PGL(3),
  \end{displaymath}
  a finite $\GL(3,\bF_2)$-quotient of a smooth projective GIT
  quotient $(\bP^2)^7//\PGL(3)$.
\end{theorem-intro}

\begin{theorem-intro}\label{thm-intro:burniat}
  For Burniat surfaces, the normalization of the 
  compactification $\oM\bur\slc$ is the quotient of a certain moduli
  space $\oM(\frac12)$ of labeled stable pairs 
  (Def.~\ref{def:Mc}) by the finite group $C_3\ltimes
  S_2^4$. The normalization map is a bijection.

  There is a diagram of moduli spaces of labeled stable pairs
  (see Section~\ref{sec:25})
  \begin{displaymath}
    \oM\tor = \oM(\tfrac13) \xleftarrow{\rho_1} \oM(\tfrac25)
    \xleftarrow{\rho_2} \oM(\tfrac12)
  \end{displaymath}
  in which $\oM\tor$ is a projective toric variety
  with $8$ isolated singularities
  corresponding to an
  explicit fan $\fF$ (Def.~\ref{def:fanF}), $\rho_1$ is the blowup
  at one smooth point, and $\rho_2$ is the blowup at six disjoint
  smooth rational curves avoiding the singular locus.
\end{theorem-intro}

\begin{remark} The moduli space of smooth surfaces of general type with fixed
numerical invariants is open in the corresponding moduli space of
stable surfaces,
but possibly not dense: here by  ``compactification'' of a class 
of surfaces we mean its closure in the moduli space of stable surfaces. 
In fact for Burniat surfaces there are  additional  components
meeting the compactification $\oM\bur\slc$ of Theorem
\ref{thm-intro:burniat} (cf. Remark \ref{rem: extra-Burniat}).
We don't know if the main component of $\oM\bur\slc$ is
normal. 

On the other hand,  in the Campedelli case we conjecture that the compactification $\oM\cam\slc$ of Theorem \ref{thm-intro:campedelli} is a connected component of the moduli space of stable surfaces.
\end{remark}

The boundary  $\oM\cam\slc\setminus M\cam$ of $M\cam$ is the union of  two irreducible divisors  and the boundary of $\oM\bur\slc$ is the union of  eight  irreducible divisors. Note that, contrary to the case of curves,  it is not always the case that the complement $\oM\slc \setminus M$ of a component $M$ of the moduli space of surfaces inside its closure in the moduli of stable surfaces is a divisor (see for instance \cite{Franciosi2017K21}, Table 1).

Here all the surfaces corresponding to boundary points are $\bZ^k_2$-covers, induced by the bicanonical system, of Gorenstein surfaces with Cartier total branch divisor, so  by the Hurwitz formula they  have Cartier index 1 or 2. This  is also 
somewhat unexpected, since although the index is known to be bounded for any fixed value of $K^2$ (\cite{alexeev1994boundedness-and-ksp-2}) already for $K^2=1$ there are examples of stable surfaces  with Cartier index 15 (\cite{coughlan2023Tdivisors}) and 21 (\cite{gallardo2022unimodal}).

\bigskip

The first version \cite{alexeev2009explicit-compactifications} of this
paper was written in 2009, and by some accounts it served as an
introduction to the subject to many students. We haven't finished
it until now, however, for two reasons:

The main reason was that it used the moduli space of stable pairs
$(X,\sum b_iB_i)$ with coefficients $b_i=\frac12$, which did not
really exist at the time. The moduli with fixed coefficients
$b_i\le\frac12$ present problems on the level of the definition of
families, as the divisors $B_i$ may form non-flat families.  Various 
 solutions, none completely satisfactory, were proposed: working
with subschemes $B_i\subset X$ instead of divisors; working with
finite maps $B_i\to X$; restricting to seminormal reduced bases, etc.
Recently, a good solution involving the notion of  K-flatness was
proposed, and a complete theory has been firmly put in place in
\cite{kollar2023families-of-varieties}.

The second reason was that the initial computation used an ad hoc
generalization of the theory of weighted hyperplane arrangements
\cite{alexeev2015moduli-weighted}, which itself was not fully worked
out at the time.

In the present version we give new, easier proofs and provide much
sharper, very explicit descriptions of the compactified moduli
spaces. We then sketch the original proofs which use the moduli of
stable weighted hyperplane arrangements.

Results of the present paper were used in
\cite{alexeev2023secondary-burniat} to compute the stable surface 
compactifications of secondary Burniat surfaces and in
\cite{alexeev2023kappa-classes} to investigate kappa classes on KSBA
spaces. 

\smallskip

We work over $\bC$ since the general results of
\cite{kollar2023families-of-varieties} about the existence of the
stable pair compactifications are known only over $\bC$. But in fact
most of the constructions work, or can be modified to work, over any
field $k$ of characteristic different from~$2$.
The source of this paper on arXiv includes a sagemath \cite{sagemath} file
verifying computations with fans and polytopes. 

\begin{acknowledgements}
  The first author was partially supported by the NSF under
  DMS-2201222. The second author is a member of GNSAGA of INdAM.
  We thank J\'anos Koll\'ar for helpful comments.
\end{acknowledgements}

\section{Preliminaries}
\label{sec:premilimaries}

\subsection{Compact moduli of stable surfaces}
\label{subsec:compact-moduli}

We briefly recall the main definitions and the existence theorem,
referring the reader to \cite{kollar2023families-of-varieties} for
more details. 

We say that a variety has double crossings if
every point is either smooth or has a neighborhood formally
isomorphic to $xy=0$. It is {\em deminormal}
if it is $S_2$ and has double crossings outside a closed subset of
codimension $\ge2$.

Let $X$ be a variety, let $B_j$, $j=1,\dots n$, be effective Weil
divisors on $X$, possibly reducible and with components in common,
and let $b_j$ be rational numbers
with $0<b_j\le 1$. Set $B=\sum_j b_jB_j$.

\begin{definition}\label{defn:lc}
  Assume that $X$ is a normal variety. Then $X$ has a canonical
  Weil divisor $K_X$ defined up to linear equivalence. The pair
  $(X,B)$ is called \emph{log canonical} (lc) if
  \begin{enumerate}
  \item $K_X+B$ is $\mathbb Q$-Cartier, i.e. some positive multiple is
    a Cartier divisor, and
  \item every prime divisor $D$ of $X$ has multiplicity $\le 1$ in $B$ and
    for every proper birational morphism $h\colon X'\to X$ with normal
    $X'$, in the natural formula
    \begin{math}
      K_{X'} + h\inv _*B = h^*(K_X+B) + \sum a_i E_i
    \end{math}
    one has $a_i\ge -1$. Here, $E_i$ are the irreducible exceptional
    divisors of $\pi$ and the pullback $h^*$ is defined by extending
    $\mathbb Q$-linearly the pullback on Cartier divisors. $h_*\inv B$
    is the strict preimage of $B$.
  \end{enumerate}

  The pair $(X,B)$ is called \emph{Kawamata log terminal} (klt) if all
  $a_i>-1$ and $\mult_D B:=\sum b_j\mult_D(B_j)<1$.
\end{definition}

\begin{remark}\label{rem: lc-for-lines}
If  $(X, \sum b_i B_i)$ is a pair such that $X$ is a smooth surface and  the $B_i$ are smooth prime divisors that intersect transversally, then Definition \ref{defn:lc} can be made more  explicit. Since the blow-up of $X$ at the points where at least 3 of the $B_i$ meet is a log resolution of $(X,B)$,  \cite[Cor.~2.23]{kollar2013singularities}  applies and $(X,\sum b_i B_i))$ is  log canonical  iff the following conditions hold:
\begin{enumerate}
\item every prime divisor  appears in  $\sum b_iB_i$ with multiplicity $\le 1$, 
\item ever point  has multiplicity $\le 2$ for  $\sum b_iB_i$.
\end{enumerate}
\end{remark}
\begin{definition}
  A pair $(X,B)$ is called \emph{semi log canonical} (slc) if 
  \begin{enumerate}
  \item $X$ is deminormal,
  \item no divisor $B_j$ contains any component of
    the double locus of $X$,
  \item some multiple of the Weil $\mathbb Q$-divisor $K_X+B$, well
    defined thanks to the previous condition, is Cartier, and
  \item denoting by $\nu\colon X^{\nu}\to X$ the normalization, the
    pair $(X^{\nu},\ \text{(double locus)} + \nu_*\inv B )$ is log
    canonical.
  \end{enumerate}
\end{definition}

\begin{definition}
  A pair $(X,B)$ is a KSBA stable pair, or simply a \emph{stable pair}, if
  \begin{enumerate}
  \item $(X,B)$ has slc singularities, and
  \item $K_X + B$ is ample.
  \end{enumerate}
\end{definition}

A family of stable pairs over a normal scheme $S$ is a flat morphism $f\colon \cX\to S$ together
with Weil divisors $\cB_i\subset\cX$ such that $K_\cX + \sum_ib_i\cB$ is a
relative ample $\bQ$-Cartier divisor and every geometric fiber is a
stable pair. For general (e.g.,  not reduced or reduced but not
seminormal) base schemes   the   definition of  families of divisors $\cB_i$ is much more 
delicate, and we refer to \cite[Ch.~7]{kollar2023families-of-varieties} for
more details. The families we construct will be over reduced normal
bases.

\begin{theorem}[\cite{kollar2023families-of-varieties}, Thm.~8.1]
  \label{thm:ksba-moduli}
  For fixed $(b_1,\dotsc,b_n)$ and fixed  $(K_X+B)^{\dim X}$, there exists a
  coarse moduli space of stable pairs, and it is projective. 
\end{theorem}

\subsection{Abelian covers}
\label{subsec:covers}

A $G$-cover is a finite morphism $X\to Y$ of varieties which is
the quotient map for a generically faithful action of a finite 

group~$G$. This means that for every component $Y_i$ of $Y$ the
$G$-action on the restricted cover $X\times_Y Y_i\to Y_i$ is faithful.
We will restrict ourselves to the case when $G$ is abelian.

When $Y$ is smooth and $X$ is normal, the theory of abelian covers was
described in \cite{pardini1991abelian-covers}. In
\cite{alexeev2012non-normal-abelian} we extended it to the case needed
for this paper: when one or both of $X$ and $Y$ are non-normal and
deminormal. We briefly review this theory. 

The $G$-action on $X$ with $X/G=Y$ is equivalent to a decomposition:
\begin{displaymath}\label{eq:decomposition}
  \pi_*\cO_X=\bigoplus_{\chi\in G^*}\cF_{\chi},\qquad \cF_0=\cO_Y
\end{displaymath}
where $G$ acts on $\cF_{\chi}$ via the character $\chi$; each sheaf $\cF_{\chi}$ is generically locally free of
rank~$1$. The variety $X$ is $S_2$ iff each sheaf $\cF_\chi$ is
$S_2$.

Now assume that $Y$ is smooth and that $X$ is $S_2$. Then 
$\cF_\chi=L\inv_\chi$ for some invertible sheaves $L_\chi$ on $Y$. The
$\cO_Y$-algebra structure on $\pi_*\cO_X$ is given by  global sections
$s_{\chi,\chi'}$ of $L_\chi\otimes L_{\chi'}\otimes
L\inv_{\chi+\chi'}$. No section $s_{\chi,\chi'}$ is identically zero, since otherwise $X$ would be non reduced, so $D_{\chi,\chi'} =
(s_{\chi,\chi'})$ is an effective divisor  and $L_\chi\otimes L_{\chi'} \simeq
L_{\chi+\chi'}(D_{\chi,\chi'})$. 

For the rest of the section we restrict to the case $G=\bZ_2^k$, which
is especially simple and is enough for our applications. When writing
down formulas for this case it is convenient to write $G$ additively
and identify the group of characters $G^*=\Hom(G,\bG_m)$ with the dual
vector space $(\bZ_2^k)^{\vee}$, so that the natural pairing takes
place in $\bZ_2$. 

Then one of the main results of
\cite{pardini1991abelian-covers} is that if $X$ is normal then there exist unique effective
divisors $D_g$ labeled by the nonzero elements $0\ne g\in G$ such that
\begin{math}
  D_{\chi,\chi'}=\sum_{g:\ \chi(g)=\chi'(g)=1} D_g.
\end{math}
The support of the divisor $D_g$ equals the image of the divisorial
part of the set of points $x\in X$ fixed by the automorphism $g$.  and
$X\to Y$ is \'etale outside of $\cup D_g$.  We will call
$\sum D_g=D\tot$, or simply $D$, the \emph{total branch divisor}. 

We remark (cf. \cite[\S~1]{pardini1991abelian-covers}) that for more general abelian groups $G$ the divisors $D_{H,\psi}$ are labeled by the cyclic subgroups $H\subseteq G$ with a choice of a generator $\psi$ of $H^*$. When $G=\bZ_2^k$, the pairs $(H,\psi)$ are obviously in a bijection with the nonzero elements of $G$.
The above description extends to the case $Y$  smooth  and $X$
deminormal, see \cite[Cor.~1.10]{alexeev2012non-normal-abelian}.
\medskip

\begin{definition}\label{def: fundrel}
The data of the invertible sheaves $L_\chi$ for all $\chi\in G^*$ and
effective divisors $D_g$ for $0\ne g\in G$ are called the
\emph{building data} of the cover; they satisfy the identities
\begin{equation}\label{eq:fundrel}
  \text{(written additively)} \quad
  L_{\chi} + L_{\chi'} \simeq L_{\chi+\chi'} + 
  \sum_{g:\ \chi(g)=\chi'(g)=1} D_g,
\end{equation}
 called the \emph{fundamental relations}. Note that for $\chi=\chi'$ one gets the relation $2L_{\chi} \simeq \sum_{g:\ \chi(g)=1} D_g$.
The divisors $D_g$ are called the {\em branch data}: if $\Pic(Y)$ has no 2-torsion the branch data suffice to determine the cover.
\end{definition}

Vice versa, given  building data satisfying the fundamental relations, if
$H^0(\cO_Y^*) = \bC^*$ (for example if $Y$ is proper and connected)
then there exists a unique cover $\pi\colon X\to~Y$ with these building
data. Without assuming $H^0(\cO_Y^*) = \bC^*$, the cover exists and is
unique \'etale-locally.

\begin{example}[$\bZ_2^2$-covers] \label{ex:burniat} Consider the group 
  $\bZ_2^2$. Denote the nonzero elements of $\bZ_2^2$ by $R,G,B$ for the
  colors red, green and blue, and use the same letters $R,G,B$ to
  denote the corresponding divisors. Then the fundamental relations
  imply that there exist  line bundles $L_1,L_2,L_3$ such that: 
  \begin{equation} \label{eq: RGB} 
  2L_1 = G+B,\quad
  2L_2=B+R, \quad 2L_3=R+G.
  \end{equation}
   Vice versa, assume that one has three
  effective divisors $R,G,B$ such that the divisors $G+B$, $B+R$, $R+G$ are
  $2$-divisible in $\Pic Y$. If $\Pic(Y)$ has no torsion the solutions $L_1,L_2, L_3$ of \eqref{eq: RGB} are unique and  $R,G,B$
  together  with $L_1,L_2,L_3$  are the building data of a $\bZ_2^2$-cover.
\end{example}

The following is easy:

\begin{lemma}\label{lem:connected}
  Suppose that $Y$ is smooth and connected and $\Pic(Y)$ has no 2-torsion. Then $X$ is connected iff 
  $G$ is generated by the elements  $g$ with $D_g\ne0$. 
\end{lemma}
\begin{proof}
Let $G_0<G$ be the subgroup generated by the $g$ such that $D_g\ne 0$ and set $X_0:=X/G_0$. Then the induced cover $X_0\to Y$ is an \'etale $G/G_0$-cover and therefore it is trivial since $\Pic(Y)$ has no 2-torsion. So if $X$ is connected then $G=G_0$. 

Conversely, if $G=G_0$ then by the fundamental relations for every $\chi\ne 0$ we have $2L_{\chi}>0$, hence $h^0(L_{\chi}\inv)=0$. So $h^0(\OO_X)=1$ by the projection formula, and thus $X$ is connected. 
\end{proof}

The extension of the theory to the case $X$, $Y$ demi-normal goes as follows (cf. \cite{alexeev2012non-normal-abelian}). If  $Y$ is singular but normal one uses $S_2$-fication: if
$i\colon U\to Y$ is the nonsingular locus, then $\codim(Y\setminus U)\ge2$,
and for any $G$-cover $\pi\colon X\to Y$ with $X$ demi-normal the restriction $\pi\inv(U)=:V\to U$ is a demi-normal 
$G$-cover with  smooth base, so  $\pi_*\cO_V=\oplus_{\chi\in G^*} L_\chi\inv$ for suitable line bundles $L_{\chi}$  on $U$. Since   $X$ is $S_2$ by assumption,  $\pi_*\cO_X=\OO_Y\oplus\cF_\chi$ with $\cF_\chi = i_*L\inv_\chi$. The
pushforwards of invertible sheaves from $U$ to $X$ are  divisorial
sheaves corresponding to Weil divisors modulo linear
equivalence. Thus, all the same statements about the building data and
fundamental relations hold, with the Weil divisors $D_g$ and
$L_\chi$ taken in the class group $\Cl(Y)$.

By computing in local coordinates (see \cite[\S
3]{pardini1991abelian-covers}), with the above reduction from the
case of a normal $Y$ to the case of a smooth $Y$, 
one has for $G=\bZ_2^k$: 
\begin{lemma}\label{lem:cover-normalization}
  Suppose that $Y$ is normal and $X$ is deminormal. Then:
  \begin{enumerate}
  \item the components of $D\tot=\sum_{g\in G} D_g$ have 
  multiplicities $\le2$;
  \item  $X$ is normal iff
  $D\tot$ is  reduced;
 \item  the normalization of~$X$ is a
  $G$-cover  $X'\to Y$ with  branch data $D'_h$,  defined as follows:
   a prime divisor
   $E$ of $Y$ is a component of $D'_h$ iff $h=\sum_g\left(mult_{D_g} E\right)g$. In particular, $E$ is not a component of $D':=\sum_hD'_h$ iff  $sum_g\left(mult_{D_g} E\right)g=0$.
  \end{enumerate}
  \end{lemma}

\medskip Finally, the case when both $X$ and $Y$ are deminormal 
is treated by \cite{alexeev2012non-normal-abelian}, Theorems~1.13 and
1.17. The main result is that every $G$-cover $\pi\colon X\to Y$ is
obtained from a $G$-cover $\pi'\colon X'\to \wY$ of the normalization
by a gluing construction.

\smallskip

Finally we recall from \cite{alexeev2012non-normal-abelian} a Hurwitz type formula for  the canonical class of $X$ in terms of the canonical class of $Y$ and of the branch data for a cover $\pi\colon X\to Y$ of deminormal varieties. Let $E$ be a prime divisor of $Y$ and set:
\begin{itemize}
\item $a_E=0$ if $\pi$ is generically \'etale over $E$ or if $E$ is contained in the double locus of $Y$,
\item  $a_E=1$ if $Y$ is generically smooth along $E$ but $X$ is singular along $\pi\inv(E)$,  
\item $a_E=\frac12$ otherwise. 
\end{itemize}The divisor $D\hur=\sum_Fa_FF$ is called the Hurwitz
divisor of $\pi$. Note that if $Y$ is normal  $D\hur=\frac12 D\tot$
and this equality holds more generally if the irreducible components
of $Y$ are smooth in codimension 1.
We  have the following
(\cite[Lem.~2.3, Prop.~2.5]{alexeev2012non-normal-abelian}):

\begin{lemma}
  \label{lem:slc-cover}
  Let $\pi\colon X\to Y$ be a $\bZ_2^k$-cover of deminormal 
  varieties. Then:
  \begin{enumerate}
  \item $2K_X=\pi^*( 2K_Y+2D\hur)$ in $\Cl(X)$;
  \item  $K_X$ is $\bQ$-Cartier iff so is $K_Y+D\hur$, and
  $X$ is slc iff so is the pair $(Y,D\hur)$.
  \end{enumerate}
\end{lemma}

\section{Campedelli surfaces with $\pi_1(X)=\bZ_2^3$}
\label{sec:Campedelli}

\subsection{Definitions}

We will work with canonical models of surfaces of general type. Thus,
our normal surfaces of general type will have canonical (i.e. Du Val)
singularities and ample canonical class.  As usual in surface theory, we write $p_g(X)=h^2(\OO_X)$ and $q(X)=h^1(\OO_X)$. 

The term \emph{(numerical) Campedelli surface} normally refers to a
surface of general type with $K_X^2=2$ and $p_g=q=0$. The first
examples of such surfaces were constructed by Campedelli
\cite{campedelli1932sopra-alcuni} in 1932.

The Campedelli surfaces with fundamental group of order 8 are usually
described as free quotients of the intersection of 4 quadrics in
$\bP^6$ by a group of order 8 (cf.  \cite{miyaoka1977on-numerical} for the
case $\pi_1=\Z_2^3$ and \cite{lopes2009campedelli-surfaces} for the general
case).  When $G=\bZ_2^3$, the quadrics can be taken to be diagonal and
it is easy to check that the bicanonical system gives a $\bZ_2^3$-cover
of $\pp^2$ branched on 7 lines. These are the surfaces we
consider. Over $\bC$ at least, they form a connected component in the
moduli space of surfaces of general type with canonical singularities.

\begin{definition}
  For brevity, a \emph{Campedelli surface} in this paper will denote a
  $\bZ_2^3$-cover $\pi\colon X\to\bP^2$ whose building data are 7 lines
  $D_g$ ($g\in\bZ_2^3\setminus 0$) for which the cover has canonical
  singularities.

  This means that either the lines are in general position (and then
  the cover is smooth) or three distinct lines $D_{g_1}$, $D_{g_2}$, $D_{g_3}$ intersect at a point
  and the three elements $g_1,g_2,g_3$ generate $G$ (in which case the
  cover has an   $A_1$ singularity, see \cite[Table~1]{alexeev2012non-normal-abelian}).

  We will denote the moduli space of Campedelli surfaces with
  canonical singularities by $\Mcamp$, and the open subset of smooth
  surfaces by $\Mcamp^0$.
\end{definition}

 The fundamental
relations~\eqref{eq:fundrel} have the solution
$L_{\chi}=\cO_{\bP^2}(2)$ for every $0\ne \chi\in G^*$. Thus,
$h^i(\cO_X)=\sum h^i(L_{\chi}\inv)=0$ for $i=1,2$.
By  Lemma~\ref{lem:slc-cover} one has $K_X = \pi^*(K_{\bP^2}+\frac12\sum D_g)= \pi^*(\frac12 h)$, $h$ being the class of a line in $\pp^2$, so
indeed $K_X^2 = 8\cdot(\frac12)^2=2$. Using the standard projection formulas for abelian covers, one can show  $|2K_X|=\pi^*|\OO_{\pp^2}(1)|$, so the covering map $\pi$ coincides with the bicanonical map. In particular, every automorphism of $X$ descends to an automorphism of $\pp^2$ that permutes the 7 branch lines of $\pi$. So for a general choice of these lines we have $\Aut(X)=G$.

Let $\M^0(3,7)$ denote the moduli space of arrangements of 7 lines in
$\bP^2$ in general position, i.e., such that no three of them are concurrent: it is the free $\PGL(3)$-quotient of an open
subset of $(\bP^2{}^\vee)^7$.
The moduli space $\Mcamp^0$ of smooth
Campedelli surfaces is obtained by dividing $\M^0(3,7)$ by the choice of a
basis in $\bZ_2^3$. Thus, the coarse moduli space is $\Mcamp^0 = \M^0(3,7)
/ \GL(3,\bF_2)$.

\subsection{Compact moduli spaces}\label{sec:compact-campedelli}
To describe the compactified moduli space, we need to understand two
things:
\begin{enumerate}
\item what are the degenerations of arrangements of 7 lines, and
\item what happens to the abelian covers.
\end{enumerate}

\begin{theorem}
  Let $F = (\bP^2)^7//\PGL(3)$
  be the GIT quotient for the ``democratic'' polarization
  $(1,\dots,1)$. Then there exists a family
  $(\cY, \sum_{i=1}^7 \frac12 \cB_i)\to F$ whose fibers are log
  canonical hyperplane arrangements
  $(\bP^2, \sum_{i=1}^7 \frac12 B_i)$. The GIT quotient is also the
  geometric quotient for the set of stable points, the group action is
  free, and the quotient is smooth and projective.
\end{theorem}
\begin{proof}
  The space $(\bP^2)^7$ parameterizes the set of ordered $7$-tuples of
  hyperplanes in the dual projective plane, and
  $\PGL(3) = \Aut \bP^2$.
  By \cite[Prop.~4.3]{mumford1994geometric-invariant} (see also
  \cite[Thm.~II.2.1]{dolgachev1988point-sets}) the pair
  $(\bP^2; B_1,\dotsc, B_7)$ is GIT stable (resp. semistable) iff the
  following two conditions hold:
  \begin{enumerate}
  \item the number of lines coinciding with a given line is $<\frac73$
    (resp. $\le\frac73$),
  \item the number of lines passing through a point is $<\frac{14}3$
    (resp. $\le\frac{14}3$).
  \end{enumerate}
  On the other hand, by Remark \ref{rem: lc-for-lines} the pair $(\bP^2,\sum\frac12
  B_i)$ is  log canonical iff: 
  \begin{enumerate}
  \item[(1$'$)] the number of lines that coincide is $\le2$,
  \item[(2$'$)] the number of lines passing through a point is $\le4$.
  \end{enumerate}
  Since $\lfloor\frac73\rfloor=2$ and $\lfloor\frac{14}3\rfloor=4$,
  these three pairs of conditions are equivalent, so the sets of stable,
  semistable and log canonical pairs are all the same. The GIT quotient 
  of the semistable locus is projective. On the other hand,
  the quotient of the stable locus is a geometric quotient and
  the stabilizers are finite. By enumeration, one checks that every
  stable configuration contains $4$ lines in general position. Thus,
  the stabilizers are trivial, the group action is free, and the
  quotient is smooth.

  The universal family $\cY$ is the quotient of the universal family
  of line arrangements over the semistable = stable locus $\big(
  (\bP^2)^7 \big)^s$ by $\PGL(3)$.
\end{proof}

\subsection{Proof of Theorem~\ref{thm-intro:campedelli}} 
\label{sec:proof-campedelli}

\begin{proof}
  For any point in $F=(\bP^2)^7//\PGL(3)$, choosing a sufficiently
  small open neighborhood $U\subset F$, one can identify
  $\cY\times_Y U$ with $\bP^2\times U$, and the sheaves $L_\chi$ and
  $\cO_\cY(D_g)$, $g\in\bZ_2^3\setminus 0$, are pullbacks from $\bP^2$.
  Then
  \begin{displaymath}
    \cX:=\Spec_\cY \left(\oplus_{\chi\in G^*} L_{\chi}\inv \right)    
  \end{displaymath}
  gives a family of semi log canonical surfaces. For each surface $X$
  in this family, $\cO(2K_X)$ is an ample invertible sheaf defining a
  $\bZ_2^3$-cover $X\to\bP^2$ (see section \ref{subsec:covers}, and in
  particular Lemma \ref{lem:slc-cover}), and the local deformations of
  $X$ are given by deforming the $7$ lines $D_g$. The lines $D_g$ are
  labeled by $0\ne g\in\bZ_2^3$, and two choices differ by a choice of
  a basis, i.e. by an element of $G=\GL(3,\bF_2)$.  Thus, there is a
  bijective map from $\oM\cam=F/G$ to coarse moduli space
  $\oM\cam\slc$. We claim that this bijection is an isomorphism to an
  irreducible component of $\oM\cam\slc$.

  The deformations of $X$ with $G$-action are equivalent to the
  deformations of the pair of $\bP^2$ with $7$ lines. So the space of
  $G$-deformations is smooth. The space of $G$-deformations is
  immersed into the space of all deformations, as the the tangent
  space to the former is the linear subspace of the tangent space to
  the latter on which $G$ acts trivially.
\end{proof}

\begin{remark}\label{rem:gerbe}
  For the ``labeled'' stable Campedelli surfaces, with the branch
  divisors labeled by $0\ne g\in\bZ_2^3$, the moduli stack is a gerbe
  over $F$ banded by $\bZ_2^3$. It would have been the quotient stack
  $[F:\bZ_2^3]$ if the sheaves $L_\chi$ had global square roots over
  $F$. Similarly, for the unlabeled stable Campedelli surfaces, the
  moduli stack is a gerbe over the stack $[F:\GL(3,\bF_2)]$ banded by
  $\bZ_2^3$.
  We thank Angelo Vistoli for explaining this point to us.
\end{remark}

\subsection{Degenerate Campedelli surfaces and their singularities}

The singularities occurring on the degenerate Campedelli surfaces were
considered in detail in Tables 1,2,3 of  \cite{alexeev2012non-normal-abelian}. Enumerating the
possibilities for the lines $D_g$, $g\in\bZ_2^3\setminus 0$ gives the
following:

\begin{lemma}\label{lem:sing_camp}
  In the notation of \cite{alexeev2012non-normal-abelian}:
  \begin{enumerate} 
  \item 
  the  singularities occurring on degenerate Campedelli
  surfaces are  $3.1$, $3.3$, $4.3$, $4.4$, $2'.1$, $3'.1$, $3'.4$,
  $4'.5$, $4'.6$, $4'.7$, $4''.4$, 
  $4''.5$. 
\item $\oM\cam$ contains two boundary divisors, one consisting of
  surfaces with an $A_1$ singularity (type $3.1$) and the other one
  consisting of surfaces with two $\frac{1}{4}(1,1)$ singularities
  (type $3.3$).
  \end{enumerate}
\end{lemma}
  Here, the notation $k.n$ means that $k$ lines pass through a common
  point on the base surface $\bP^2$; $n$ is a case number from
  \cite{alexeev2012non-normal-abelian}. 
  Similarly, $k'.n$ means that
  two of them coincide to form a double line, and $4''.n$ means that
  there are two pairs of double lines. The integer $n$ refers to the
  $n$-th case in the list of possible relations between the $k$
  lines.

\begin{proof}
  The proof of (1) is a direct enumeration of cases.
  Part (2) is a consequence of the fact that three branch lines
  $D_{g_1}$, $D_{g_2}$, $D_{g_3}$ going through the same point is the
  only codimension one degeneration. When $g_1,g_2,g_3$ are linearly
  independent in $\bZ_2^3$, the cover has an $A_1$ singularity. When
  $g_1+g_2+g_3=0$, it has two  $\frac14(1,1)$ singularities.
\end{proof}

\begin{remark}\label{rem:campedelli-all-cases}
  It is a straightforward but tedious exercise to list the boundary
  data of higher codimensions.
  Indeed, over an infinite field of $\chr k\ne2$, modulo $S_7$ there are $36$
  configurations of $7$ lines in $\bP^2$ such that $\le 2$ lines
  coincide at a time and $\le 4$ lines pass through a common
  point. Modulo our relabeling group $\GL(3,\bF_2)\subset S_7$ there
  are $175$ orbits.
  It is not very practical to list them all here.

  The cases of codimension $2$ are as follows:
  \begin{enumerate}
  \item The lines $D_1,D_2,D_3$ pass through a common point, and the
    lines $D_4,D_5,D_6$ pass through a common point. There are two
    cases: $g_1+g_2+g_3=0$ and $\{g_4,g_5,g_6\}$ is a basis of
    $\bZ_2^3$, or both $\{g_1,g_2,g_3\}$ and $\{g_4,g_5,g_6\}$ are
    bases. 
  \item The lines $D_1,D_2,D_5$ pass through a common point,
    and the lines $D_3,D_4,D_5$ pass through a common point. 
    There are four orbits depending on the triples $ijk \in
    \{125, 345, 567\}$ for which $g_i+g_j+g_k=0$:
    all; $125$; $567$; none. In the last two cases the singularities
    are the same: two $A_1$.
  \item Two lines $D_1=D_2$ coincide. There is only one case,
    with the singularities $2'.1$, $3'.1$, $3'.4$ of Table 2 in
  \cite{alexeev2012non-normal-abelian}.
  \item Four lines pass through a common point. Again, there are two
    cases mod $\GL(3,\bF_2)$.
    This gives the normal, log canonical but not log terminal
    singularities, cases $4.3$ 
    and $4.4$ of Table 2 in \cite{alexeev2012non-normal-abelian}.
  \end{enumerate}
 \end{remark}

\begin{lemma} \label{lem:camp_components}
  A degenerate Campedelli surface may have 1, 2, or 4 (but not 8)
  irreducible components.
\end{lemma}
\begin{proof}
  By Lemmas~\ref{lem:cover-normalization}, \ref{lem:connected}, a
  $\bZ_2^3$-cover is not normal iff the branch divisor $\sum_gD_g$ is
  not reduced and  the normalization is also a $\bZ_2^3$-cover
  branched on a divisor contained in $\sum_gD_g$. In addition, the slc
  condition implies that no three of the $D_g$ can coincide. If
  $D_g=D_{g'}=h$, the line $h$ occurs in the branch locus of the
  normalization with label $g+g'\ne 0$. So the normalization is
  branched on at least four lines and the case of 8 components never
  occurs.
  
    Up to the action of $\GL(3,\bF_2)$, the case of four irreducible components occurs when
  $D_{100}=D_{011}$, $D_{010}=D_{101}$, $D_{001}=D_{110}$, where we
  use the natural labels for the nonzero elements of $\bZ_2^3$. In
  this case the normalization has 4 components, each of them a
  double cover of $\bP^2$ ramified in 4 lines corresponding to
  $g=111$. Each component is a del Pezzo surface of degree 2 with  
 six  $A_1$ singularities. It is easy to see that up to the action of
  $\GL(3,\bF_2)$ this is the only case with 4 irreducible components.

If we split one of the double lines then the cover has 2 components. Each of
them is a del Pezzo of degree~1 with six $A_1$ singularities.
\end{proof}

\section{Burniat surfaces with $K_X^2=6$}
\label{sec:Burniat}
In this paper, we consider only Burniat surfaces with $K_X^2=6$. These
are the so called ``primary'' Burniat surfaces. There exist also
``secondary'' and ``tertiary'' Burniat surfaces with $K_X^2=5$ and
$4$. They were considered in \cite{hu14compactifications-of-moduli} and \cite{alexeev2023secondary-burniat}. 

\subsection{Set-up and notation}

\begin{figure}[htbp!] 
  \centering
  \includegraphics{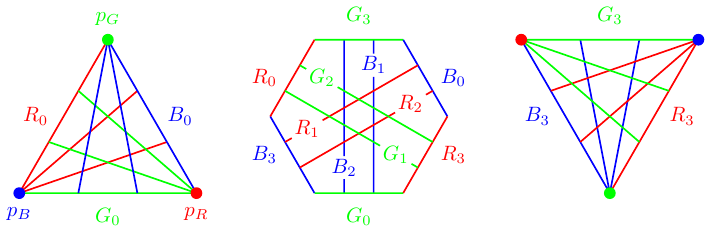}      
  \caption{Burniat arrangements on $\bP^2$ and $\Sigma=\Bl_3\bP^2$}
  \label{fig-burniat-config}
\end{figure}

Consider the arrangement of 9 lines on $\bP^2$ shown in the first panel
in Fig.~\ref{fig-burniat-config}.
Using the RGB color scheme ($R$ for red, $G$ for green, $B$ for blue;
in the black and white version of this paper they are shown as solid,
dashed and dotted) 
we denote the sides of the
triangle $R_0,G_0,B_0$ and the vertices $p_R,p_G,p_B$. The point $p_R$
is the point of intersection of $G_0$ and $B_0$, etc. There are
additional lines $R_1,R_2$ through $p_B$, lines $G_1,G_2$ through
$p_R$, and lines $B_1,B_2$ through $p_G$. We assume that the lines are
in general position otherwise.

We note the cyclic RGB symmetry of this picture. 

Now blow up the points $p_R,p_G,p_B$ and denote the resulting
exceptional divisors on the surface $\Sigma=\Bl_3\bP^2$ by
$R_3,G_3,B_3$.  Note that the arrangement on $\Sigma$ can be presented
as the blowup of $\bP^2$ in a different way by contracting
$R_0,G_0,B_0$. The two line arrangements differ by a Cremona
transformation.

\begin{definition}
  Set $R=\sum_{i=0}^3 R_i$, $G=\sum_{i=0}^3 G_i$,
  $B=\sum_{i=0}^3 B_i$, corresponding to the 3 nonzero elements of
  $\bZ_2^2$.  The divisors $R+G$, $G+B$ and $B+R$ are $2$-divisible in
  $\Pic(\Sigma)$, the fundamental relations \eqref{eq:fundrel} have a
  unique solution (cf. Example \ref{ex:burniat}), and there is a
  (unique) $\bZ_2^2$-cover 
  $\pi\colon X\to \Sigma$ with branch data $R,G,B$. The surface $X$ is
  called a Burniat surface.
\end{definition}

When the lines are chosen generically,  so that on $\Sigma$ only two
divisors  at a time intersect (and they belong to different elements of
$G=\bZ_2^2$, which is always true for Burniat arrangements), the Galois
cover is smooth.
 In the notation of section \ref{subsec:covers}, we have  $D\hur=\frac 12 D\tot=\frac12(R+G+B)$, 
so $K_{\Sigma}+ D\hur = -\frac12 K_{\Sigma}$ is ample, and by Lemma \ref{lem:slc-cover}
 $K_X =
\pi^*( K_{\Sigma}+ D\hur )$, and so $K_X^2 = 4\cdot \frac64 =6$.
To compute $p_g(X)=q(X)=0$ one solves equations
\eqref{eq:fundrel} for $L_{\chi}$  and uses the projection formula
$\pi_*\cO_X=\cO_\Sigma\oplus_{G^*\setminus 0}L_{\chi}^{-1}$.

By \cite{lopes2001connected-component} Burniat surfaces form a connected
component $\Mbur$ in the moduli space of canonical surfaces of general
type. The dimension of $\Mbur$ is 4. The map to $\Sigma$ is the
bicanonical map, so it is intrinsic to the surface $X$. 

\begin{definition}\label{def:relabeling-group}
  We define the \emph{relabeling group} for the tuple
  $(\Sigma, R_i,G_i,B_i)$ to be
  $\Gamma = C_3\ltimes S_2^4 \subset S_{12}$, a group of order $48$,
  acting as follows:
  \begin{enumerate}
  \item $S_2^3$ acts by exchanging $R_1\leftrightarrow R_2$,
    $G_1\leftrightarrow G_2$, $B_1\leftrightarrow B_2$ independently.
  \item The remaining copy of $S_2$ acts by an involution exchanging the six
    $(-1)$-curves, the sides of the hexagon $R_0\leftrightarrow R_3$,
    $G_0\leftrightarrow G_3$, $B_0\leftrightarrow B_3$ at the same
    time. 
  \item The cyclic group $C_3$ rotates the RGB colors.
  \end{enumerate}
  This is the subgroup of the automorphism group of the dual graph of
  $\cup R_i\cup G_i\cup B_i$ that preserves the partition
  into three parts, arbitrarily colored. 
\end{definition}

\begin{lemma}\label{lem:label-group}
  The automorphism group of a smooth Burniat surface $X$ with the labeled curves
  $R_i,G_i,B_i$ is the covering group $\bZ_2^2$ of $\pi\colon
  X\to\Sigma$.   
    The smooth Burniat surfaces $X$, $X'$ defined by two Burniat
  configurations $(\Sigma, R_i,G_i,B_i)$ and
  $(\Sigma, R'_i,G'_i,B'_i)$ are isomorphic iff there exists an
  automorphism $\alpha\colon\Sigma\to\Sigma$ sending the $12$ $RGB$
  curves to the $12$ $R'G'B'$ curves and such that the induced
  permutation of the $12$ labels is an element of $\Gamma$.
\end{lemma}
\begin{proof}
  The automorphism group of $\Sigma$ preserving the $6$ boundary
  curves is the torus $\bG_m^2$ and the subgroup of this torus acting
  trivially on two non-parallel boundary curves is trivial. The
  automorphism group of $\bP^1$ fixing $3$ points is trivial. For a
  smooth Burniat surface there are even $4$ distinct points on each
  boundary curve. 
  Thus an automorphism of a labeled Burniat surface
  acts trivially on $\Sigma$, so it is an element of the covering group.

  Since the bicanonical map $\varphi_{|2K_X|}\colon X\to\Sigma$ is
  intrinsic, any isomorphism $X\to X'$ induces an isomorphism
  $\alpha\colon\Sigma\to\Sigma$ permuting the $12$ branch curves.
   It
  must send the $4$ curves in each branch divisor $D_g$ to curves in
  some $D_{g'}$, so it should permute the three colors. With the
  colors fixed, obviously $\alpha(R_1)=R_1$ or $R_2$, etc. and
  $\alpha(R_0)=R_0$ or $R_3$. And if $\alpha(R_0)=R_3$ then
  $\alpha(G_0)=G_3$ and $\alpha(B_0)=B_3$ since the boundary hexagon
  is distinguished. So $\alpha\in \Gamma$.
\end{proof}

We note the following obvious observation (cf. Remark \ref{rem: lc-for-lines}).  To simplify the notation,
we use $R_i$, $G_i$, $B_i$ to denote both the lines in $\bP^2$ and
their strict transforms in $\Sigma$.
\begin{lemma}\label{lem:lc-conditions}
  Let $\big(\bP^2,\sum_{i=0}^2 (r_iR_i+g_iG_i+b_iB_i)\big)$ and
  $\big(\Sigma,\sum_{i=0}^3 (r_iR_i+g_iG_i+b_iB_i)\big)$ be pairs as
  in Fig.~\ref{fig-burniat-config}, $f\colon\Sigma\to\bP^2$ be the
  blowup, and suppose that
  \begin{displaymath}
    f^*\left(K_{\bP^2}+\sum_{i=0}^2 (r_iR_i+g_iG_i+b_iB_i)\right) =
    K_\Sigma + \sum_{i=0}^3 (r_iR_i+g_iG_i+b_iB_i).
  \end{displaymath}
  Then
  \begin{displaymath}
    r_3=g_0+g_1+g_2+b_0-1, \ g_3=b_0+b_1+b_2+r_0-1, \ b_3=r_0+r_1+r_2+g_0-1
  \end{displaymath}
  and the first pair is log canonical iff so is the second one.
\end{lemma}

\subsection{Variation of weights}
\label{sec:burniat-variation}

The curves on $\Sigma$ are split into two groups: \emph{boundary} and
\emph{interior}:
\begin{displaymath}
  D\bry = R_0+R_3+G_0+G_3+B_0+B_3, \quad
  D\inr = R_1+R_2+G_1+G_2+B_1+B_2.
\end{displaymath}
We have $D\bry\equiv -K_{\Sigma}$ and $D\inr\equiv -2K_{\Sigma}$. Thus, the
$\bQ$-divisor $K_{\Sigma} + \frac12 D\bry + cD\inr$ is ample for any
$c>\frac14$.

\begin{definition}\label{def:Mc}
  For $c>\frac14$, we denote by $\oM(c)$ the normalization of the
  compactification for the moduli space of pairs
  $(\Sigma,\frac12 D\bry + cD\inr)$, which exists by
  Theorem~\ref{thm:ksba-moduli}, with the labeled and ordered curves
  $R_i,G_i,B_i$.
\end{definition}

We are ultimately interested in $\oM(\frac12)$ but we will proceed in
stages. By Remark \ref{rem: lc-for-lines} the non log canonical singularities of the pair $(\Sigma,\frac12
D\bry + c D\inr)$ occur in the interior $\Sigma\setminus D\bry$ when
\begin{enumerate}
\item $c > \frac13$ and $6$ of the interior lines meet at a point of $\Sigma\setminus D\bry$, or
\item $c > \frac25$ and $5$ of the interior lines meet at a point of $\Sigma\setminus D\bry$.
\end{enumerate}
Since two lines of the same color meet at an interior point iff they coincide, condition (1) above is equivalent to $R_1=R_2$, $G_1=G_2$, $B_1=B_2$ and $R_1,G_1,B_1$ meet at an interior point,  and condition (2) can be rephrased in a similar way. 

For $\frac14 < c \le \frac13$ the moduli space $\oM(c)$ is the same. 
For these weights the singularities are possibly not log
canonical only if some of the curves $R_i,G_i,B_i$ with $i=1,2$ go to
the boundary. These degenerations are purely toric.

\subsection{The toric setup}
\label{sec:toric-setup}

We fix the torus embedding $\bG_m^2 \into \Sigma$. The coordinates on
$\bG_m^2$ can be chosen symmetrically to be $x,y,z$ with $xyz=1$. Then
the divisor $D\inr$ on $\bG_m^2$ is given by the equation
\begin{equation}\label{eq:f}
  F = (x+r_1)(x+r_2)(y+g_1)(y+g_2)(z+b_1)(z+b_2).
\end{equation}
The Newton polytope of $F$ is a side-$2$ hexagon, and $F$ defines a
section of $\cO_\Sigma(2)$.

We choose the orientation in such a way that $r_i\to 0$ for $i=1,2$
means $R_i\to R_0+G_3$, $g_i\to0$ means $G_i\to G_0+B_3$,
$b_i\to 0$ means $B_i\to B_0+R_3$. 

There is a natural action of $\bG_m^6$ on the equation $F$, rescaling
$r_i,g_i,b_i$, and the torus action on $\Sigma$ gives an
embedding $\bG_m^2\to \bG_m^6$. It gives an exact sequence of tori
\begin{displaymath}
  1 \to N_\Sigma\otimes\bG_m \to N_\cY\otimes\bG_m \to N\otimes \bG_m\to 1,
  \quad \text{where}
\end{displaymath}
\begin{enumerate}
\item
  $N_\cY=\bZ^6=\{v=(\rho_1,\rho_2,\gamma_1,\gamma_2,\beta_1,\beta_2)\}$,
\item the sublattice $N_\Sigma\simeq\bZ^2 \subset N_\cY$ is the set of
  vectors $(\rho,\rho,\gamma,\gamma,\beta,\beta)$ with
  $\rho+\gamma+\beta=0$,
\item We identify
  $N$ with the subset of $\bZ^4$ of
  quadruples $(\delta,\bar\rho,\bar\gamma,\bar\beta)$ such that
  $\delta+\bar\rho+\bar\gamma+\bar\beta \equiv 0\pmod 2$ and the map
  $N_\cY\to N$ with
  \begin{displaymath}
    v \mapsto (\rho_1+\rho_2+\gamma_1+\gamma_2+\beta_1+\beta_2,
    \rho_1-\rho_2,\gamma_1-\gamma_2,\beta_1-\beta_2).
  \end{displaymath}
\end{enumerate}

We will define two toric varieties and a toric morphism $\cY\to
\oM(\frac12)$ by explicit fans $\fF_\cY$ in $N_\cY$ and $\fF$ in
$N$ and a map of fans $\fF_\cY\to\fF$.

\begin{lemma}\label{lem:relabeling-group}
The relabeling group $\Gamma=C_3\ltimes S_2^4$ of Definition~\ref{def:relabeling-group} acts
as follows:
\begin{enumerate}
\item $S_2^3$ acts on $N_\cY$ by switching the order in each of the pairs
  $(\rho_1,\rho_2)$, $(\gamma_1,\gamma_2)$, $(\beta_1,\beta_2)$,
  and on $N$ by sending $(\delta,\bar\rho,\bar\gamma,\bar\beta)$ to
  $(\delta,\pm\bar\rho,\pm\bar\gamma,\pm\bar\beta)$.
\item another $S_2$ acts by sending
   $v\to -v$ and
  $(\delta,\bar\rho,\bar\gamma,\bar\beta)\to(-\delta,-\bar\rho,-\bar\gamma,-\bar\beta)$. 
\item $C_3$ acts by cyclically permuting the groups $(\rho_1,\rho_2)$,
  $(\gamma_1,\gamma_2)$, $(\beta_1,\beta_2)$, and by permuting the
  three coordinates 
  $\bar\rho,\bar\gamma,\bar\beta$ in $(\delta, \bar\rho,\bar\gamma,\bar\beta)$.
\end{enumerate}
\end{lemma}

\subsection{Minimal (codimension~$1$) toric degenerations}
\label{sec:mindegs}

Consider a DVR $R$ with quotient field $K$ and a valuation
$\nu\colon K^*\to\bZ$ with a generator $t$ of the maximal ideal
$\mathfrak m=(t) \subset R$, $\nu(t)=1$.
Without loss of generality we can as well
take $R=\bC[t]_{(t)}$ or $\bC[[t]]$. The computations are the same but
notation is easier. Instead of a family over $\Spec K\subset\Spec R$
one can equivalently work with a family over the germ $(\Delta,0)$ of
a smooth curve.

\begin{notations}\label{not:1-param-family}
Any one parameter degeneration of
the Burniat configuration  is described by a $6$-tuple
$(r_1,r_2,g_1,g_2,b_1,b_2)\in (K^*)^6$.
We can write $r_i = r_i'\cdot t^{\rho_i}$, $g_i = g_i'\cdot
t^{\gamma_i}$, $b_i = b_i'\cdot t^{\beta_i}$
with $\rho_i=\nu(r_i)$, $\gamma_i=\nu(g_i)$, $\beta_i=\nu(b_i)$ and
with $r_i',g'_i,b'_i$ invertible in $R$.
Thus, any one-parameter degeneration defines a vector
$(\rho_1,\rho_2,\gamma_1,\gamma_2,\beta_1,\beta_2)\in N_\cY$ and
its image $(\delta,\bar\rho,\bar\gamma,\bar\beta)\in N$.
Denote the residues of $r'_i,b'_i,g'_i$ in $R/\mathfrak m=\bC$
by $\barr_i,\barg_i,\barb_i\in\bC^*$.
\end{notations}

\begin{figure}[htbp!] 
  \centering
  \includegraphics[width=360pt]{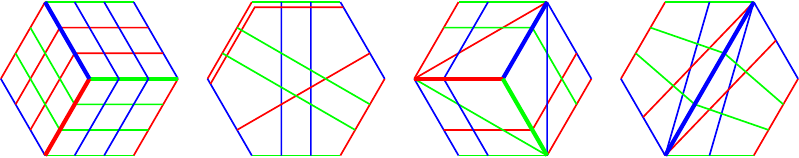}
  \caption{Minimal toric degenerations of types A, B, C, D.}
  \label{fig-mindegs}
\end{figure}

We begin by describing the four ``minimal'' one-parameter degenerations shown  in Fig.~\ref{fig-mindegs},  which we call A, B, C and
D.  In the figures  the thicker lines denote the curves in the double locus of the limit surface, and, as usual in toric geometry, a triangle represents $\pp^2$,  a rhombus represents  $\pp^1\times\pp^1$,  a trapezoid represents $\mathbb F_1$ and a hexagon the del Pezzo surface of degree 6. So in  Fig.~\ref{fig-mindegs} B is irreducible, A and C consist of three copies of $\pp^1\times \pp^1$ and D is the union of two copies of $\mathbb F_1$. 

Each of these four degenerations  produces a $3$-dimensional family of limit pairs.  The
cases are ordered by the slope
$\mu = \frac{|\delta|}{|\bar\rho|+|\bar\gamma|+|\bar\beta|}$.  The orbits
of these one-parameter degenerations under the relabeling group
$\Gamma$ will then define the rays of the fans $\fF_\cY$ and
$\fF$.

\ifshort
\else
\fi

\subsubsection*{Case A}
$(1,1,0,0,0,0), (0,0,1,1,0,0), (0,0,0,0,1,1) \mapsto (2,0,0,0)$.

\begin{figure}[htbp!] 
  \centering
  \includegraphics[width=317pt]{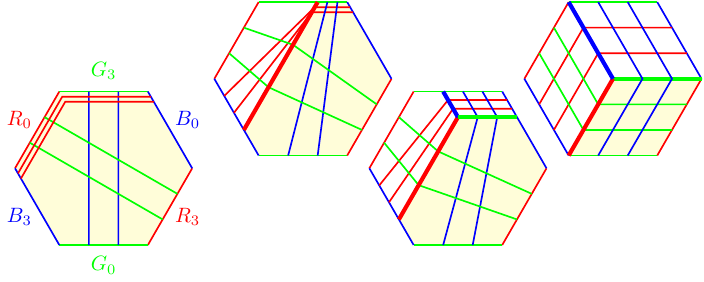} 
  \caption{Degeneration for case A}
  \label{fig-A}
\end{figure}

In the limit, $R_i\to R_0+G_3$ for $i=1,2$. The limits of the divisors
$G_i,B_i$, $i=1,2$ remain in the interior of $\Sigma$.  We have a
constant family of varieties $Y=\Sigma\times\Spec R$ and $12$ divisors
$R_i,G_i,B_i$ on it.

Blow up the line $R_0$ in the central fiber $Y_0$. Then the central fiber
becomes $\Bl_3\bP^2\cup \bF_1$.  Blowing up the strict preimage of $G_3$ 
changes $\bF_1$ into $\Bl_2\bP^2$ and inserts
$\bF_0=\bP^1\times\bP^1$.
To make such computations, we use the well-known \emph{triple point
  formula:} Let  $Y=\cup Y_j$ be the central fiber in a smooth
one-parameter family, and assume that $Y$ is reduced and simple normal
crossing. Let $C$ be the intersection $Y_1\cap Y_2$, suppose it is
smooth. Then
\begin{displaymath}
  (C|_{Y_1})^2+(C|_{Y_2})^2+ 
  \text{(the number of the triple points of $Y$ contained in $C$)}
  = 0.
\end{displaymath}

For the central fiber, the divisor $K_Y+\frac12 D\tot$ restricted to an
irreducible component $Y_j$ is $K_{Y_j}+ \frac12 D\tot|_{Y_j}+ \text{(the double
locus)}|_{Y_j}$. 
The curves $R_i,G_i,B_i$ appear in the last sum with
coefficient $\frac12$, and the curves in the double locus with
coefficient~$1$.

A simple computation shows that after the last step on the central
fiber $K_Y+\frac12 D\tot$ is big, nef and zero on 3 curves. The 3-fold
pair $(\cY,\frac12 \cD\tot)$ is simple normal crossing.
The Basepoint-Free Theorem
\cite[Thm.~3.24]{kollar1998birational-geometry} immediately implies
that a big positive multiple $N(K_{\cY}+\frac12\cD\tot)$ gives a
birational morphism contracting the three zero curves.

After the contraction the new central fiber is a union of three
$\bP^1\times\bP^1$ together with 8 curves on each. The equations of
$D\inr$ on these components are
\begin{displaymath}
  \prod_{i=1,2}(y+\barg_i)(z+\barb_i), \quad
  \prod_{i=1,2}(z+\barb_i)(x+\barr_i), \quad
  \prod_{i=1,2}(x+\barr_i)(y+\barg_i). 
\end{displaymath}
There is a $3$-dimensional family of
such pairs, parameterized by $(\bC^*)^3$, and all of them appear as
limits. For example, we can take $\barr_2=\barg_2=\barb_2=1$ and vary
$\barr_1$, $\barg_1$,~$\barb_1$. 
A natural compactification of $(\bC^*)^3$ in $\oM(\frac12)$
is $(\bP^1)^3$, sending either $R_1$ or $R_2$ to the
boundary of the hexagon, and similarly for $G_i$ and $B_i$. 

\begin{remark}\label{rem:colors-double-locus}
  The labeling of the double locus is done in such a way that the
  fundamental relations \eqref{eq:fundrel} of Definition \ref{def: fundrel} hold on each irreducible
  component, so one has a well defined $\bZ_2^2$-cover. Namely, on
  each irreducible component the divisors $R+G$, $G+B$ and $B+R$ must
  be divisible by $2$.  Formally, the new color is obtained by
  multiplying in the group $\bZ_2^2$ the colors intersecting at the
  blown-up locus. We will justify this choice in
  Lemma~\ref{lem:double-curve-colors}.
\end{remark}

\subsubsection*{Case B}
$(1,0,0,0,0,0) \mapsto (1,1,0,0)$. 
In this case the limit of the curve $R_1$ coincides with $R_0+G_3$.
The resulting configuration is still log canonical and there is
obviously a $3$-dimensional family of limits parameterized by $(\bC^*)^3$.

\medskip

The remaining two cases are very similar to case A. We let the
pictures do the explanations for the blowups and blowdowns.

\subsubsection*{Case C}
$(0,-1,1,0,1,0), (1,0,0,-1,1,0), (1,0,1,0,0,-1) \mapsto
(1,1,1,1)$.  
\begin{figure}[htbp!] 
  \centering
  \includegraphics[width=317pt]{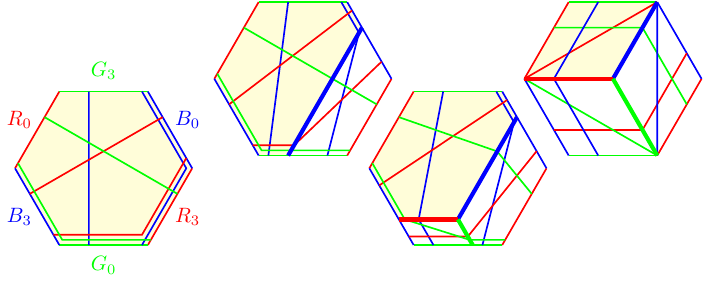} 
  \caption{Degeneration for case C}
  \label{fig-C}
\end{figure}

Again, the central fiber $Y$ is a union of three
$\bP^1\times\bP^1$. The equations of $D\inr$ restricted to these irreducible
components are (remember that $xyz=1$):
\begin{displaymath}
  (x+\barr_1)(y+\barg_2)(z+\barb_2), \ \ 
  (y+\barg_1)(z+\barb_2)(x+\barr_2), \ \ 
  (z+\barb_1)(x+\barr_2)(y+\barg_2).
\end{displaymath}
For each component the moduli space is $\bC^*$, 
giving $(\bC^*)^3$ as the parameter space. Again, we can take
$\barr_2=\barg_2=\barb_2=1$ and vary $\barr_1$, $\barg_1$, $\barb_1$
to realize all of them.

\begin{figure}[htbp!] %
  \centering
  \includegraphics[width=261pt]{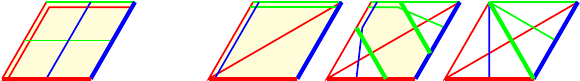} 
  \caption{Further case C degenerations over $0$ and $\infty\in\bP^1$}
  \label{fig-C+}
\end{figure}

The natural compactification of each $\bC^*$ is $\bP^1$. 
There are two stable degenerations shown in
Fig.~\ref{fig-C+}.  The degenerations of the three components are
independent since the gluings are unique. So the compactification of
this family in the stable pair moduli space is $(\bP^1)^3$.

\subsubsection*{Case D}
$(1,0,0,0,0,-1), (0,-1,0,0,1,0) \mapsto (0,1,0,1)$
\begin{figure}[htbp!]
  \centering
  \includegraphics[width=246pt]{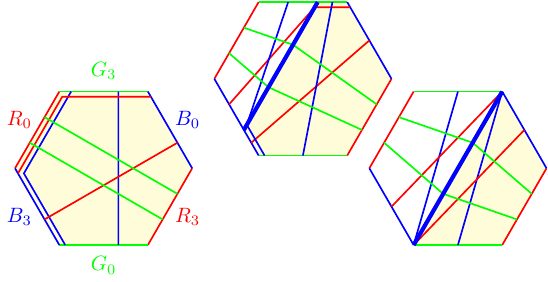} 
  \caption{Degeneration for case D}
  \label{fig-D}
\end{figure}

For each of the two irreducible components, isomorphic to $\bF_1$,
there is a $2$-dimensional moduli space, but they should reduce to the
same $4$ points in $\bP^1$ on the intersection, giving a moduli space
of dimension $2+2-1=3$ of such surfaces.

Indeed, the equations of $D\inr$ on the components are
(remember that $xyz=1$):
\begin{displaymath}
  (x+\barr_2)(z+\barb_1)(y+\barg_1)(y+\barg_2), \ \ 
  (x+\barr_1)(z+\barb_2)(y+\barg_1)(y+\barg_2).
\end{displaymath}
Taking 
$\barr_2=\barg_2=\barb_2=1$ and varying $\barr_1$, $\barg_1$,
$\barb_1$ realizes the parameter space $(\bC^*)^3$.

\begin{remark}
 It is not hard to give explicit equations of a degeneration of the del Pezzo surface  $\Sigma$ of degree 6 to the  underlying surface $Y_0$ of the limit pair for the minimal degenerations of type A,C and D (Figure \ref{fig-mindegs}). 
 
 For type A and C consider the threefold $\mathcal Y\subset \left(\pp^1\right)^3\times \bA^1_t$ defined by $x_1x_2x_3+ty_1y_2y_3=0$, where $x_i,y_i$ are  homogeneous coordinates on the $i$-th copy of $\pp^1$, $i=1,2,3$.  Denote by $Y_t$ the fiber over $t\in \bA^1$ of the projection $\mathcal Y\to \bA^1$;
 if $t\ne 0$ the surface $Y_t$ is isomorphic to $\Sigma$, while the special fiber $Y_0$ is the union of three copies of $\pp^1\times \pp^1$ intersecting along rulings as shown in Figure \ref{fig-mindegs}. The singularities of $\mathcal Y$ are three ordinary double points, one on each component of the double curve of $Y_0$.
 
 For type D consider the $\pp^2$-bundle $\pp:=\Proj\left(\OO_{\pp^1}\oplus\OO_{\pp^1}\oplus\OO_{\pp^1}(1)\right)$ with relative homogeneous coordinates $x,y,z$. Denote by $u,v$ the homogeneous coordinates of $\pp^1$ and let  $\mathcal Y \subset \pp\times \bA^1_t$ be the threefold defined by $xy+t(uvz^2)=0$. Standard computations show that for $t\ne0$  the fiber $Y_t$ over $t\in \bA^1$ of the projection $\mathcal Y\to \bA^1$ is a smooth del Pezzo surface of degree 6, hence it is isomorphic to $\Sigma$, while the special fiber $Y_0$ is a union of two copies of $\bF_1$ meeting transversally along a section of self-intersection $1$.  The singularities of $\mathcal Y$ are two ordinary double points, both on the double curve of $Y_0$.

\end{remark}

\subsection{Maximal (codimension~$4$) toric degenerations}
\label{sec:maxdegs}

\begin{figure}[h!]
  \centering
  \includegraphics[width=360pt]{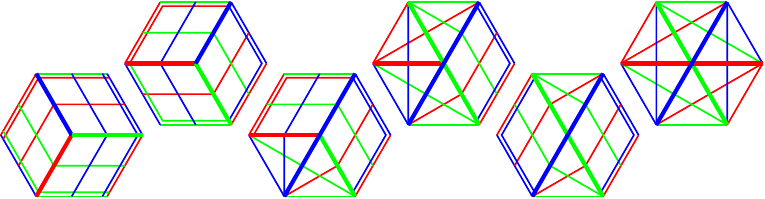}
  \caption{Maximal toric degenerations of the Burniat arrangement}
   \label{fig-maxdegs}
\end{figure}

Next, we describe $6$ maximal degenerations, obtained as combinations
of the minimal degenerations of Fig.~\ref{fig-mindegs}. They are
shown in Fig.~\ref{fig-maxdegs}. The last surface is a cone over the
cycle of rational curves. 

Each of these degenerations naturally defines a
maximal-dimensional cone in the $4$-dimensional lattice $N$, and we
list the integral generators of its rays.
\begin{enumerate}[align=left, leftmargin=5em]
\item[(AB$^3$)] $(2,0,0,0)$, $(1,1,0,0)$, $(1,0,1,0)$, $(1,0,0,1)$.
\item[(B$^3$C)] $(1,1,0,0)$, $(1,0,1,0)$, $(1,0,0,1)$, $(1,1,1,1)$.
\item[(B$^2$CD)] $(1,1,0,0)$, $(1,0,0,1)$, $(1,1,1,1)$, $(0,1,0,1)$.
\item[(BCD$^2$)] $(1,0,0,1)$, $(1,1,1,1)$, $(0,1,0,1)$, $(0,0,1,1)$.
\item[(B$^2$D$^2$)] $(1,0,0,1)$, $(-1,0,0,1)$, $(0,1,0,1)$, $(0,0,1,1)$.
\item[(C$^2$D$^3$)] $(1,1,1,1)$, $(-1,1,1,1)$, $(0,0,1,1)$, $(0,1,0,1)$,
  $(0,1,1,0)$. 
\end{enumerate}
Note that the first five cones are nonsingular, and the last one is
non-simplicial.

\subsection{The ad hoc fans $\fF$ and $\fF_\cY$}
\label{sec:adhock-fans}

We will define the toric family $\cY\to\oM\tor$ and the corresponding
map of fans $\fF_\cY\to\fF$ more intrinsically in
Section~\ref{sec:toric-family}, but that theory has a high entry
point. However, it is very easy and instructive to define the fans
$\fF$, $\fF_\cY$ directly from the minimal and maximal degenerations
already found. That is what we do here, calling them the ``ad hoc
fans''. We will check in Proposition~\ref{prop:Fs-are-correct} that
they are indeed the same as the fans of
Section~~\ref{sec:toric-family}.

\begin{definition}\label{def:fanF}
  We define $\fF$ to be the fan in $N$ with
  \begin{enumerate}
  \item The rays $\bR_{\ge0}(2,0,0,0)$, $\bR_{\ge0}(1,1,0,0)$,
    $\bR_{\ge0}(1,1,1,1)$, $\bR_{\ge0}(0,1,0,1)$ for the minimal
    degenerations 
    of Section~\ref{sec:mindegs}, and their images under the
    relabeling group $\Gamma=C_3\ltimes S_2^4$. We call them the rays
    of types A, B, C, D respectively.
  \item Maximal cones of the six types listed in
    Section~\ref{sec:maxdegs} $+$ their images under~$\Gamma$.
  \end{enumerate}
  We define $\oM\tor$ to be the toric variety with the fan $\fF$.
\end{definition}

We checked the following by hand and also in sage:

\begin{lemma}\label{lem:fanF}
  $\fF$ is a complete fan with $(1,42,210,328,160)$ cones of dimensions
  $0\text{--}4$. It is the normal
  fan of the polytope $\Pi(a,b,c)$ with  vertices
  $(a, -a + b + 4c, -a + b + 4c, -a + b + 4c)$, \ 
  $(3c, b + c, b + c, b + c)$,\  $(2c, b, b + 2c, b + 2c)$, \ 
  $(c, b + c, b + c, b + 3c)$,\  $(0, b, b, b + 4c)$,\ 
  $(0, b + 2c, b + 2c, b + 2c)$ and their images under the action of
  $\Gamma$, for any $a,b,c\in\bR$ with $b,c>0$ and
  $b+4c > a > 3c$.
\end{lemma}

Note: the fact that the fan $\fF$ is complete implies that we got all
the degenerations.

\begin{lemma}\label{lem:cones}
  The $\Gamma$-orbits of the cones in $\fF$ of dimensions $1$--$4$ are
  as follows:
\emph{
  \begin{enumerate}
  \item A($2$), B($12$), C($16$), D($12$).
  \item AB($12$), B$^2$($6$), B$^2$($24$), BC($48$), BD($24$), BD($24$), 
   CD($48$), D$^2$($24$).
  \item AB$^2$($24$), B$^3$($16$), B$^2$C($48$), B$^2$D($24$),
    B$^2$D($12$), B$^2$D($12$),
    BCD($48$), BCD($48$),
    BD$^2$($48$),
    CD$^2$($48$). 
  \item AB$^3$($16$), B$^3$C($16$), B$^2$CD($48$), BCD$^2$($48$),
    B$^2$D$^2$($24$), C$^2$D$^3$($8$). 
  \end{enumerate}
  }
\end{lemma}

\begin{corollary}\label{cor:M-sings}
  The toric variety $\oM\tor$ for the fan $\fF$ is projective. It has
  $8$ isolated singularities corresponding to the {\rm C$^2$D$^3$}
  cones. Each singularity is isomorphic to the cone over
  $\big(\bP^1\times\bP^2, \cO(1,1)\big)$ and admits two small
  resolutions.
\end{corollary}

We now define the fan $\fF_\cY$ for the presumed family
$f\colon\cY\to\oM\tor$ from general principles of toric geometry:
The rays of $\fF_\cY$ should correspond to the torus-invariant divisors
in $\cY$. The maximal cones of $\fF_\cY$ should correspond to the
torus-fixed points in $\cY$, mapped by $f$ to the torus-fixed points
in $\oM\tor$. So under the projection $N_\cY\to N$ the maximal cones
of $\fF_\cY$ must map to the maximal cones of~$\fF$.

\begin{definition}\label{def:fF_Y}
  The torus-invariant divisors in $\cY$ are of two kinds:

  \smallskip
  
  (1) Those that correspond to the irreducible components of the
  minimal degenerations. Under $f$ they map to the divisors of types
  A, B, C, D in $\oM\tor$. 
  The corresponding rays of $\fF_\cY$ are listed in
  Section~\ref{sec:mindegs}. For example, there are the A-type rays
  $\bR_{\ge0}(1,1,0,0,0,0)$, $\bR_{\ge0}(0,0,1,1,0,0)$,
  $\bR_{\ge0}(0,0,0,0,1,1)$ mapping to the A-type ray
  $\bR_{\ge0}(2,0,0,0)$ of~$\fF$.
  Thus, for each type-A ray of $\fF$ there are three A-type rays of
  $\fF_\cY$ mapping to it, for B-type there is one, for C-type there are
  three, and for D-type there are two.
  
  \smallskip

  (2) The divisors and their rays corresponding to the $6$ boundary
  curves which do not vary in the family:
  \begin{eqnarray*}
    R_0\, (-1,-1,0,0,1,1), \quad
    G_0\, (1,1,-1,-1,0,0), \quad
    B_0\, (0,0,1,1,-1,-1), \\
    R_3\, (1,1,0,0,-1,-1), \quad
    G_3\, (-1,-1,1,1,0,0), \quad
    B_3\, (0,0,-1,-1,1,1).
  \end{eqnarray*}
  These divisors dominate $\oM\tor$, and the rays map to $(0,0,0,0)\in N$.
  We will call these the divisors and rays of type $\Delta$.

  \smallskip

  We define the rays of $\fF_\cY$ to be these rays and their images
  under $\Gamma=C_3\ltimes S_2^4$, for a total of $2\cdot 3 + 12 +
  16\cdot 3 + 12\cdot 2 + 6 = 96$ rays.

  \medskip
  
  For each surface in Fig.~\ref{fig-maxdegs} there are $7$ torus-fixed
  points, corresponding to the $7$ vertices of the polytopes. For each
  of them we list the irreducible components of the minimal
  degenerations whose closures contain this point, and the boundary
  divisors containing it. The corresponding rays define a
  maximal-dimensional cone of $\fF_\cY$. We define the maximal
  cones of $\fF_\cY$ to be their images under $\Gamma$.
\end{definition}

\begin{lemma}\label{lem:F_Y-cones}
  For the maximal cones of $\fF$ of the six types listed in
  Section~\ref{sec:maxdegs}, the following are the types of maximal
  cones of $\fF_\cY$ mapping to them:
  \begin{enumerate}[align=left, leftmargin=2em]
    {\rm
\item[(AB$^3$)] A$^3$B$^3$, three A$^2$B$^3${$\Delta^2$}, three AB$^3${$\Delta^2$}.
\item[(B$^3$C)] B$^3$C$^3$, three B$^3$C$^2${$\Delta^2$}, three B$^3$C$\Delta^2$.
\item[(B$^2$CD)] B$^2$C$^3$D$^2$, two B$^2$CD$\Delta^2$, two
  B$^2$C$^2$D$\Delta^2$, B$^2$C$^2$D$^2${$\Delta^2$}, B$^2$CD$^2${$\Delta^2$}.
\item[(BCD$^2$)] BC$^3$D$^4$, BCD$^2${$\Delta^2$}, BC$^2$D$^2${$\Delta^2$}, two
  BC$^2$D$^3${$\Delta^2$}, two BCD$^3${$\Delta^2$}.
\item[(B$^2$D$^2$)] B$^2$D$^4$, two B$^2$D$^2${$\Delta^2$}, four B$^2$D$^3${$\Delta^2$}.
\item[(C$^2$D$^3$)] C$^6$D$^6$, six C$^3$D$^4${$\Delta^2$}.
  }
\end{enumerate}
\end{lemma}

One can check with sagemath that the fan $\fF_\cY$ is complete.  This
also follows by observing that for each maximal cone $\tau\in\fF$
the cones defined above cover its preimage in $N_\cY\otimes\bR$.  We
give a more direct proof in Proposition~\ref{prop:Fs-are-correct} by
showing that $\fF_\cY$ is the normal fan of a convex polytope.

\subsection{Moduli interpretation of a fiber fan}
\label{sec:fiber-fan}

In Section~\ref{sec:toric-family}
we will define the fans $\fF$, $\fF_\cY$ and the toric family
$f\colon\cY\to\oM\tor$ more intrinsically.  We will show that the ad
hoc fan $\fF$ defined in the previous section is in fact an
instance of a fiber fan, which is well known to have a moduli
interpretation.  One consequence of this fact is an explicit description of
the family of ``varying'' divisors $D\inr$, in addition to the
``fixed'' boundary divisor $D\bry$. 

\smallskip

Let $\phi\colon M_P\simeq\bZ^{n_P}\to M_Q\simeq\bZ^{n_Q}$ be an affine map and
$A\subset M_P$ a finite set. Consider two polytopes $P=\Conv A$ and
$Q=\Conv \phi(A)$, assume maximal-dimensional. Recall that a lattice
polytope $P$ defines a toric variety with an ample line bundle
as follows: $(V_P, L_P) = (\Proj S_P, \cO(1))$, where
the graded algebra $S_P$ is
\begin{displaymath}
  S_P = \oplus_{d\ge0} H^0(V_P, \cO(d)) = \oplus \bC e^{(d,m)} \quad\text{with }
  (d,m) \in \Cone(1,P)\cap\bM_P, 
\end{displaymath}
$\bM_P = \bZ\oplus M_P$.
The above data defines two projective toric varieties $T_P\acts (V_P,L_P)$
and $T_Q\acts (V_Q,L_Q)$ and a finite morphism $j\colon V_Q\to V_P$
such that $j^*(L_P)=L_Q$.

In this situation, Billera and Sturmfels
\cite{billera1992fiber-polytopes} defined the fiber polytope $\fib$ in
the lattice $\ker(M_P\to M_Q)$ as the Minkowski
integral of the fibers $P_q$, $q\in Q$. They also prove that it is a
weighted Minkowski sum of finitely many fibers, over the barycenters
of an appropriate subdivision of $Q$.  The faces of $\fib$ are in a
bijection with the coherent tilings $Q=\cup (Q_i, A_i)$,
$Q_i=\Conv\phi(A_i)$ for some subsets $A_i \subset A$. The vertices of
$\fib$ correspond to the tight tilings.

The fiber fan is the normal fan of the fiber polytope, defining a
toric variety $V_\fib$.  There exist (at least) four different
moduli interpretations of $V_\fib$:

\begin{enumerate}
\item As the ``Chow quotient'' $V_P///T_Q$
  \cite{kapranov1991quotients-toric}. This is the closure in the Chow
  variety of $V_P$ of the $T_P$-orbit of the cycle $[j(V_Q)]$. Of
  course the $T_P$-action factors through the action of the quotient torus
  $T_P/j^*(T_Q)$.

\item As a toric Hilbert scheme of $V_P$
  \cite{peeva2002toric-hilbert, haiman2004multigraded-hilbert}.

\item As a moduli space of stable toric varieties with a finite
  morphism to $V_P$ \cite{alexeev2002complete-moduli},
  \cite[Sec.~2.5]{alexeev2015moduli-weighted},
  \cite{alexeev2010complete-moduli}.

\item As the target of a morphism of  toric varieties $V_{P+\fib} \to V_\fib$.
\end{enumerate}

This paper is not the right place to discuss the common parts and the
differences between these approaches. (In fact, they are all
equivalent in our particular case, the key property being that the
sets $(1,\phi(A_i))$ span the semigroups 
$\Cone(1,Q_i)\cap\bM_Q$.)  We simply take the fourth approach:
it is the easiest and sufficient for our purposes.

The normal fan of the Minkowski sum $P+\fib$ comes with two maps to
the normal fans of $P$ and $\fib$, defining two projections
$p_1\colon V_{P+\fib}\to V_P$ and $p_2\colon V_{P+\fib}\to V_\fib$ and
a finite morphism $V_{P+\fib}\to V_P\times V_\fib$.  The second
projection gives an equidimensional family over $V_\fib$ which may
have non-reduced fibers in general.

Let $Y_t$ be a fiber in this family over a point $t\in V_\fib$. The
torus orbit $\orb(t)\subset V_\fib$ containing $t$ corresponds to a
face of the polytope $\fib$ and, by the above, to a tiling
$Q=\cup (Q_i,A_i)$.  The irreducible components $Y_i$ of $Y_{\rm red}$
are toric varieties $V_{Q_i}$ for the polytopes $Q_i$ in this tiling.

For each fiber, the first projection $Y_t\to V_P$ is a finite morphism,
which may not be an embedding in general, and the images of the
irreducible components may be non-normal.
Each fiber $Y_t$ comes with two divisors:
\begin{enumerate}
\item The ``fixed'' Weil boundary divisor $D\bry$ corresponding to the
  boundary~$\partial Q$.
\item The ``varying'' divisor $D\var$, the pullback of the Cartier
  divisor on $V_P$ that is the zero divisor of the section $\sum_{a\in A}
  e^a\in H^0(V_P,L_P)$. 
\end{enumerate}

We mention the following fact which we do not use but which may help
the reader understand some of the features of our construction: Under
some additional conditions (for example, when all the fibers are
reduced and the map $\phi|_A\colon A\to M_Q$ is injective), the fibers
$(Y,D\bry + \epsilon D\var)$ are stable pairs for $0<\epsilon\ll1$.
But we want the pair $(Y,\frac12 D\bry + (\frac14+\epsilon)D\inr)$ to
be stable instead. So what we do is a kind of a ``weighted'' version
of the ``standard'' construction.

\subsection{$\fF$ is a fiber fan, and a family of  pairs over $\oM\tor$}
\label{sec:toric-family}

We are now ready to state our construction.

\begin{figure}[htbp!]
  \centering
  \includegraphics[width=360pt]{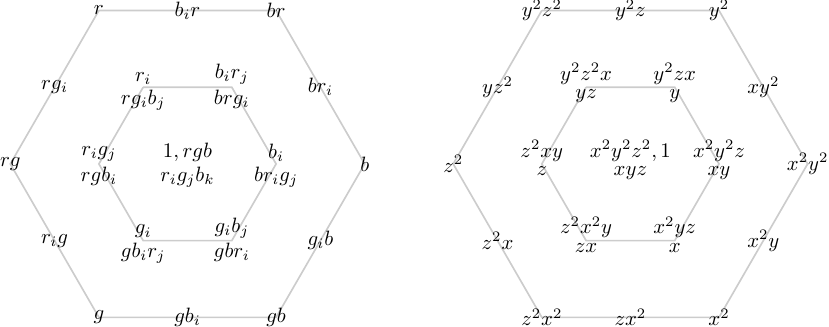}
  \caption{Monomials in $F(r_i,g_i,b_i; x,y,z)$}
  \label{fig-monomials}  
\end{figure}

\begin{definition}\label{def:rgb-projection}
  In Section~\ref{sec:toric-setup} we defined two co-character
  lattices $N_\cY=\bZ^6$ and
  $N_\Sigma=\{\rho,\gamma,\beta | \mid \rho+\gamma+\beta = 0\}$ 
  supporting the fans of  $\cY$ and $\Sigma$, with a natural inclusion
  $i\colon N_\Sigma\to N_\cY$.  The dual lattices $M_\Sigma=\bZ^6$ and
  $M_\cY=\bZ^3/\bZ(1,1,1)$ are lattices of monomials, supporting
  polytopes of projective toric varieties. We define
  an affine linear map $\phi\colon M_\cY\to M_\Sigma$ such that
  $-\phi^*\colon N_\Sigma\to N_\cY$ equals $i$. (All of our
  fans and polytopes are centrally symmetric, so $-\phi$ is chosen for
  convenience.)

  For $a=(k_1,\dotsc,k_6) \in M_\cY$, we set $\phi(a) =
  (2-k_1-k_2, 2-k_3-k_4, 2-k_5-k_6) = (\ell_1,\ell_2,\ell_3)$. We
  interpret the monomial $e^a$ as 
  \begin{math}
    r_1^{k_1}r_2^{k_2} g_1^{k_3}g_2^{k_4} b_1^{k_5}b_2^{k_6}
    \cdot x^{\ell_1}y^{\ell_2}z^{\ell_3}
  \end{math}
  and the monomial $e^{\phi(a)}$ as $x^{\ell_1}y^{\ell_2}z^{\ell_3}$
  with $xyz=1$, so that the above $rgb$-$xyz$ monomial is projected to
  its $xyz$-part.

  Let $A(F)\subset M_\cY$ be the set of the $2^6=64$ monomials
  appearing in the polynomial $F(r_i,g_i,b_i; x,y,z)$ of
  Equation~\eqref{eq:f}. Its projection $\phi \left(A(F)\right)$ is
  the set of the integral point of a hexagon of side~$2$, which we call
  $2Q$. We list the $rgb$- and $xyz$-parts of $e^a$ for $a\in A(F)$ in
  Fig.~\ref{fig-monomials}.  For example, $b_1x^2y^2z$ gives $b_1$ and
  $x^2y^2z$.  As a shortcut, we denote $r=r_1r_2$, $g=g_1g_2$,
  $b=b_1b_2$. The notation $r_i$ means that one has to repeat the
  monomial for $r_1$ and $r_2$, etc.
  
  Let $A\subset A(F)$ be the subset of $46$ points mapping to the $7$
  lattice points of the small, side-$1$ hexagon $Q$.  Let $P=\Conv A$.  By
  definition, we have $Q=\Conv\phi(A)$.
\end{definition}

\begin{proposition}\label{prop:Fs-are-correct}
  The normal fans of the polytopes $\fib$ and $P+\fib$ coincide with
  the ad hoc fans $\fF$, $\fF_\cY$, and the morphism
  $V_{P+\fib}\to V_\fib$ is identified with the toric family
  $\cY\to\oM\tor$ of Section~\ref{sec:adhock-fans}.
\end{proposition}
\begin{proof}
  We computed $\fib$ using the description given in
  \cite[Cor.~2.6]{billera1992fiber-polytopes} as the convex hull of a
  set of explicit vectors $\Phi_\Delta$, as $\Delta$ go over the
  triangulations of $Q$ with vertices in the multiset $\phi(A)$. Then
  we confirmed that $18\fib$ is a translate of the polytope
  $\Pi(6,5,1)$ defined in Lemma~\ref{lem:fanF}, whose normal fan is
  $\fF$.  (In fact, the Minkowski sum of certain $6$ fibers $P_q$ already
  gives $\fF$, but we don't need this.)

  We computed $P+\fib$ and its normal fan in sage and confirmed that
  the normal fan coincides with $\fF_\cY$. Its $160\cdot 7$ maximal
  cones are exactly the same as those described in
  Lemma~\ref{lem:F_Y-cones}.  A sage file with these computations is
  included in the source of this paper on arXiv.
\end{proof}

\begin{lemma}
  The morphism $f\colon\cY\to\oM\tor$ is flat with reduced fibers.
\end{lemma}
\begin{proof}
  From the definition, the images of maximal cones in
  $\fF_\cY$ are maximal cones in $\fF$, and it is easy to check that
  the image $p(\sigma)$ of any cone $\sigma\in\fF_\cY$ is a cone of
  $\fF$. The integral generators of $\sigma$ map to the integral generators of
  $p(\sigma)$, and it is easy to check that for any cone
  $p(\sigma)\in\fF$ the integral generators of the rays generate
  $p(\sigma)\cap N$. This implies that
  $p(\sigma\cap N_\cY)= p(\sigma)\cap N$.  By
  \cite[Thm.~2.1.4]{molcho2021universl-stacky}, which is an extension
  of Kato's flatness criterion
  \cite[Prop.~4.1]{kato1989logarithmic-structures}, this implies that
  the morphism $f$ is flat with reduced fibers.
\end{proof}

\begin{lemma}\label{lem:toric-very-ample}
  The morphism $\cY\to V_P \times \oM\tor$ is a closed embedding, and
  for every fiber $Y$ of $f\colon\cY\to\oM\tor$ the restriction of the
  line bundle $p_1^*(L_P)$, thereafter denoted by $\cO_Y(1)$, is very
  ample. The boundary divisor $D\bry\sim -K_Y$ on each fiber is given
  by a section of $\cO_Y(1)$, and so is Cartier.
\end{lemma}
\begin{proof}
  Using the torus action and the openness of very ampleness, it is
  enough to check the statement for the fibers over the torus-fixed
  points of $\oM\tor$, which are listed in
  Fig.~\ref{fig-maxdegs}. The statements follows because for each
  irreducible component $Y_i$ in these fibers the linear system
  generated by the sections $s_{a}$ corresponding to the vertices of
  the polytopes $Q_i$ is very ample. Indeed, there are only two types
  of components in Fig.~\ref{fig-maxdegs}: $(\bP^1\times\bP^1,
  \cO(1,1))$ and $(\bP^2,\cO(1))$ corresponding to a square and a
  triangle. 
  The boundary divisor $D\bry$ is given by the section corresponding
  to the central point of the hexagon $Q$. And it is easy to check that
  $K_Y+D\bry \sim 0$, same as for an ordinary toric variety.
\end{proof}

As we recalled at the end of Section~\ref{sec:fiber-fan}, the varying
divisor $D\var$ given by the ``standard'' construction is
not the internal divisor $D\inr$ that we are after.  Indeed, $D\var$
corresponds to the polytope $Q$ and is defined by a section of
$L_P=\cO(1)$, but $D\inr$ should correspond to $2Q$ and 
a section of $\cO(2)$.

\begin{definition}\label{def:Dinr}
  Let $A_0 = \{a_0\in A \mid \phi(a_0)=\text{the center of } Q\}.$ The
  $10$ monomials $e^{a_0}$, $a_0\in A_0$ are
  $\{x^2y^2z^2, r_ig_jb_k\cdot xyz, r_1r_2g_1g_2b_1b_2\}$.  We have
  $A_0\subset A \subset A(F)$.  The convex hull of $A(F)$ is the
  Newton polytope of $F$.

  It is easy to check that for any $a_0\in A_0$, one
  has $a_0+A(F) \subset A+A$ and $a_0+\Newton(F)\subset 2P$.  Thus,
  $e^{a_0}F$ defines a global section of $\cO_{V_P}(2)$.
  For any variety $Y_t\to V_P$ in our family, the restrictions of the
  $10$ sections $e^{a_0}F$ are sections of $\cO_{Y_t}(2)$ which are
  multiples of each other.
   If one of them is not identically zero
  on each irreducible component of $Y_t$ then it defines a Cartier
  divisor $D\inr$ on $Y_t$ and also a relative Cartier divisor
  $\cD\inr$ on the family $\cY\times_{\oM\tor} U\to U$ for some open
  neighborhood $U\ni t$.

  \end{definition}

\begin{proposition}\label{prop:family}
  For any variety $Y$ in the family $\cY\to\oM\tor$, the Cartier
  divisor $D\inr$ of Definition~\ref{def:Dinr} is well defined and
  thus there is a relative Cartier divisor $\cD\inr$. With this divisor
  included, the pairs $(Y,D\bry,D\inr)$ in the family $\cY\to\oM\tor$
  are the same as the pairs described in Sections~\ref{sec:mindegs},
  \ref{sec:maxdegs}.
\end{proposition}
\begin{proof}
  We will freely use
  Notation~\ref{not:1-param-family} and the notation introduced in
  Section~\ref{sec:fiber-fan} and Definition~\ref{def:rgb-projection}.
  Let $Y_t:=Y=\cup Y_i$ be a variety in our fiber fan family, and let
  $Q=\cup(Q_i,A_i)$ be the tiling corresponding to the orbit
  $\orb(t)\subset V_\fib$, as explained in
  Section~\ref{sec:fiber-fan}, so that $Y_i\simeq V_{Q_i}$.  We claim
  that there exists $a_0\in A_0$ such that the section
  $e^{(1,a_0)}\in H^0(V_P,\cO(1))$ is not identically zero on each
  irreducible component $Y_i$.  This is equivalent to the statement
  that for any face of the fiber polytope and the corresponding tiling
  $Q=\cup (Q_i,A_i)$, each set $\phi(A_i)$ contains the central point
  of $Q$.  It is sufficient to check this for maximal degenerations,
  shown in Fig.~\ref{fig-maxdegs}, which follows because the central
  point $0\in Q$ is a vertex of each~$Q_i$ and one always has
  $\operatorname{Vert} Q_i\subset \phi(A_i)$.

  Pick $a_0\in A_0$ as above.  The height function $h|_{A_0}$ achieves
  a minimum at this point. $F$ contains $e^{(1,a_0)}$, so the
  restriction of $e^{(1,a_0)}F$ to an 
  irreducible component $Y_i$ contains the monomial corresponding to
  central point of $Q$ with a coefficient which is a nontrivial
  polynomial in $\bar r'_j,\bar g'_j,\bar b'_j$
  (see Notation~\ref{not:1-param-family}). On the other hand, it is a
  product of linear terms $x+\bar r'_j$, $y+\bar g'_j$, $z+\bar b'_j$,
  $x$, $y$, $z$. So it is not identically zero for any
  $\bar r'_j,\bar g'_j,\bar b'_j$. Below, we compute these
  restrictions explicitly.

  Now consider a one-parameter family of pairs
  $(\cY^K,\cD^K\bry, \cD^K\inr)\to\Spec K$. We identify
  $(\cY^K,\cD^K\bry)$ with the constant family of toric varieties
  $\Sigma^K:=\Sigma\times\Spec K$. We get a map
  $p_1\colon \Sigma^K\to V_P^K$ and the pullback
  $p_1^*\colon H^0(V_P^K,\cO(1)) \to H^0(\Sigma^K, \cO(1))$:
  \begin{displaymath}
    \oplus_{a\in A} K e^{(1,a)} \to
    \oplus_{m\in Q\cap M_\Sigma} K e^{(1,m)},\quad
    e^{(1, k_i)} \to
    r_1^{k_1}r_2^{k_2} g_1^{k_3}g_2^{k_4} b_1^{k_5}b_2^{k_6}
    \cdot x^{\ell_1}y^{\ell_2}z^{\ell_3}
  \end{displaymath}
  with $r_i,g_i,b_i\in K$ and $x^{\ell_1}y^{\ell_2}z^{\ell_3} =
  e^{(1,\phi(a))}$. 
  We call the valuation
  \begin{displaymath}
    \nu(r_1^{k_1}r_2^{k_2} g_1^{k_3}g_2^{k_4} b_1^{k_5}b_2^{k_6}) =
    k_1\rho_1+k_2\rho_2+k_3\gamma_1+k_4\gamma_2+k_5\beta_1+k_6\beta_6
    \in \bZ
  \end{displaymath}
  \emph{the height} $h(a)$ of the monomial $a=(k_i)$.

  Thus, every point $a\in A$ defines a point $(a, h(a))$ in
  $M_\Sigma \oplus\bZ$. The convex hull of the set of these points for $a\in
  A$ projects to the hexagon $Q\subset M_\Sigma\otimes\bR$, and
  projections of the facets of the lower envelope are the polytopes
  $Q_i$ corresponding to the limit surface by the fiber fan
  construction.
  The completed family $\cY\to\Spec R$ is $\Proj S$, where $S$ is the
  $R$-subalgebra of $S_Q\otimes K$ generated by the monomials
  $t^{h(a)}e^{(1,\phi(a))}$, and the central fiber is $\Proj S/tS$.
  
  The limit of divisor $\cD\inr^K$ is obtained by restricting the
  section $e^{a_0}F$ to $\cY$ and reducing it mod $t$. This is
  computed as follows. We extend the tiling $Q=\cup Q_i$ by dilation
  to the tiling $2Q=\cup 2Q_i$. The polytopes $2Q_i$ are
  projections of the facets of the lower envelope defined above, extended by
  linearity away from the center.

  For each monomial $a=(k_i)\in F(A)$ we compute its height $h(a)$. If the
  point $(a,h(a))$ lies above the extended lower envelope, then
  $e^{(1,a)}$ reduces to $0$ mod $t$. If it lies on the extended lower
  envelope then it reduces to the monomial
  \begin{math}
    \bar r_1'{}^{k_1} \bar r_2'{}^{k_2} \bar g_1'{}^{k_3} \bar g_2'{}^{k_4}
    \bar b_1'{}^{k_5} \bar b_2'{}^{k_6} \cdot x^{\ell_1}y^{\ell_2}z^{\ell_3},
  \end{math}
  with $\bar r'_i,\bar g'_i,\bar b'_i$ introduced in 
  Notation~\ref{not:1-param-family}.

  \smallskip

  We now make concrete computations for the minimal degenerations.
  In each case we lift a vector $(\delta, \bar\rho, \bar\gamma,
  \bar\beta) \in N$ to a vector
  $(\rho_1,\rho_2,\gamma_1,\gamma_2,\beta_1,\beta_2)\in N_\cY$. 
  Two lifts differ by an element of $N_\Sigma$, i.e. a linear
  function on $M_\Sigma$, and lead to the same result. 
  
  Concretely, we choose the lifts of the rays to be the first vectors
  in the cases A, B, C, D of Section~\ref{sec:mindegs}:
  $(1,1,0,0,0,0)$, $(1,0,0,0,0,0)$, $(0,-1,1,0,1,0)$,
  $(1,0,0,0,0,-1)$. The answers for the minimal heights $h(a)$ over
  each point $\phi(a)\in 2Q$ are given in Fig.~\ref{fig-heights}, for
  all the monomials in $F$.
  \begin{figure}[htbp!]
    \centering
    \includegraphics[width=360pt]{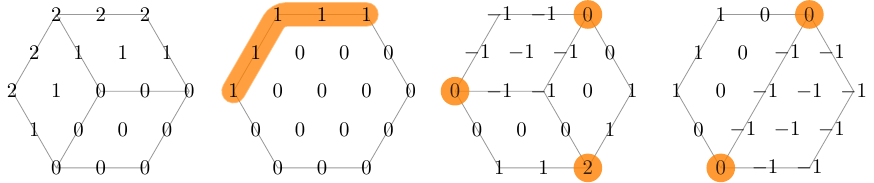}
    \caption{The minimal heights for the rays A, B, C, D}
    \label{fig-heights}
  \end{figure}
  The minimal heights for $a\in A$ are in the central $7$ points, for
  the hexagon $Q$ of side~$1$.  The picture shows the dilated domains
  of linearity, the polytopes $2Q_i$.  For the highlighted places,
  \emph{all} monomials of $F$ mapping to it vanish. For the
  non-highlighted places, some monomials survive.

  In more detail, for each polytope $2Q_i$ we compute the sum of all
  monomials in $F$ whose heights lie on the lower envelope, and get
  the following answers:
  \begin{enumerate}
  \item[(A)] $x^2 \cdot (y+g_1)(y+g_2)(z+b_1)(z+b_2) $
    and two more by symmetry.
  \item[(B)] $x(x+r_2)(y+g_1)(y+g_2)(z+b_1)(z+b_2)$
  \item[(C)] $(x+r_1)(y+g_2)(z+b_2) \cdot yzr_2$
    and two more by symmetry.
  \item[(D)] $(x+r_2)(y+g_1)(y+g_2)(z+b_1) \cdot xb_2$
    and one more by symmetry.
  \end{enumerate}
  After plugging in $r_i=r'_i t^{\rho_i}$,
  $g_i=g'_i t^{\gamma_i}$, $b_i=b'_i t^{\beta_i}$, and
  recalling the generators of the algebra $S$ above, we see that the
  equations of $D\inr$ restricted to the irreducible components of
  $\Proj S/tS$ are
  \begin{enumerate}
  \item[(A)] $(y+\bar g'_1)(y+\bar g'_2)(z+\bar b'_1)(z+\bar b'_2) $
    and two more by symmetry.
  \item[(B)] $x(x+\bar r'_2)(y+\bar g'_1)(y+\bar g'_2)(z+\bar b'_1)(z+\bar b'_2)$
  \item[(C)] $(x+\bar r'_1)(y+\bar g'_2)(z+\bar b'_2)$
    and two more by symmetry.
  \item[(D)] $(x+\bar r'_2)(y+\bar g'_1)(y+\bar g'_2)(z+\bar b'_1)$
    and one more by symmetry.
  \end{enumerate}
  So they are exactly the same as in Section~\ref{sec:mindegs}.
  Thus, the limits computed
  by the fiber polytope technology, with our ``weighted'' twist, are
  exactly the same as those that we obtained in 
  Section~\ref{sec:mindegs} by doing the Minimal Model Program steps.
  This proves the statement for the rays of $\fF$.
  
  For a vector $h$ in the interior of arbitrary cone $\tau\in\fF$ with
  rays $r_k$ and $a\in A(F)$, the height $h(a)$ is a positive combination of
  the heights $h_k(a)$. This implies
  that the tiling for $\tau$ is the intersection of the tilings for
  $r_k$, and the set of the nonvanishing monomials is the intersection
  of such sets for $r_k$. One checks that the restriction of $F$ to
  each component $Y_i$ is a product of several linear terms
  $x+\bar r'_j$, $y+\bar g'_j$, $z+\bar b'_j$, $x$, $y$, $z$
  and that the degenerations given
  by the fiber polytope technology are the same as in
  Sections~\ref{sec:mindegs}, \ref{sec:maxdegs}.
\end{proof}

\begin{warning}
  The interior divisor $\cD\inr$ is the zero set of $F$, so it is a
  Cartier divisor.  Over the open subset $(\bC^*)^4 \subset\oM\tor$ it
  splits into the sum of six divisors $R_1+R_2+G_1+G_2+B_1+B_2$. Since
  $\cD\inr$ is the closure of this set, it naturally splits into the
  sum of six Weil divisors. However, one can see from the pictures for
  the minimal degenerations C and D that sometimes these Weil divisors
  are not $\bQ$-Cartier.
\end{warning}

\subsection{Nontoric degenerations}
\label{sec:nontoric}

\subsubsection*{Case E}
In the toric family $\cY\to\oM\tor$ this degeneration occurs over the
unit point $1\in \bG_m^4$. For the pair $(\Sigma,\frac12 D\bry + c
D\inr)$ this configuration is log canonical only when $c\le\frac13$,
and is not log canonical when $c>\frac13$. Consider the latter case. 

\begin{figure}[htbp!]
  \centering
  \includegraphics[width=360pt]{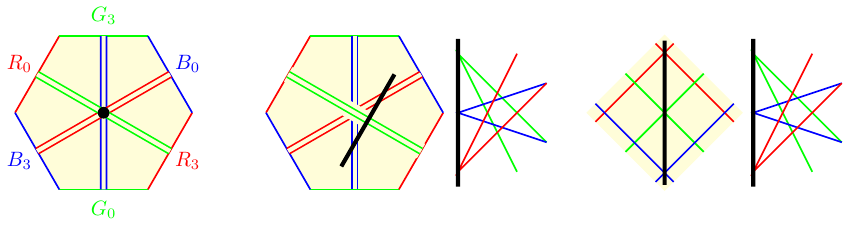} 
  \caption{Degeneration for case E}
  \label{fig-E}
\end{figure}

For a one-parameter degeneration $Y\to (\Delta,0)$ with $\Delta$
immersed into $\bG_m^4$, the resolution is obtained by blowing up a
point in the central fiber. Then the central fiber becomes
$Y_0= \Bl_p\Sigma \cup \bP^2$ with a configuration of lines on $\bP^2$
which is shown in the central panel of Fig.~\ref{fig-E}.
The intersection number for
the strict preimages of the curves $R_i,G_i,B_i$, $i=1,2$ with
$K_Y+\frac12 D\bry + c D\inr$ on the $3$-fold $Y$ is $1-2c$. It
follows that the divisor $K_{Y_0}+\frac12 D\bry + c D\inr$ is ample
for $\frac13 < c <\frac12$.

For $c=\frac12$ it is big, nef and
semiample by the standard Basepoint-Free Theorem
\cite[Thm.~3.24]{kollar1998birational-geometry}. It then
defines a contraction of three curves, the strict transforms on the blow up of $\Sigma$ of $R_1=R_2$, $B_1=B_2$, $G_1=G_2$. So on the resulting stable
model the central fiber is $Y_0 = (\bP^1\times\bP^1)\cup\bP^2$, as
shown in the right most panel  of Fig.~\ref{fig-E}. 

Now consider the configuration of lines on the $\bP^2$ in the right
panel with the condition that not all of the lines
$R_1,R_2,G_1,G_2,B_1,B_2$ pass through the same point. Such a
configuration is described by $6$ lines with the equations $x+r_1z$,
$x+r_2z$, $y+g_1z$, $y+g_2z$, $-x-y+b_1z$, $-x-y+b_2z$ modulo a free
action by the matrix group
\begin{displaymath}
  G = \left\{
  \begin{pmatrix}
    1 & 0 & e_1\\
    0 & 1 & e_2\\
    0 & 0 & 1
  \end{pmatrix}
  \cdot
  \begin{pmatrix}
    1 & 0 & 0\\
    0 & 1 & 0\\
    0 & 0 & \lambda
  \end{pmatrix}
  \qquad \text{with } e_1,e_2\in \bC, \ \lambda\in \bC^*.
  \right\}
\end{displaymath}
It is easy to see that the quotient is $(\bA^6/\bA^2)/\bG_m \simeq \bP^3$
and that it is naturally identified with the projectivization of the
tangent space of $1\in\bG_m^4$. Indeed, an easy computation shows that
the limits of the families with the immersed bases
$(\Delta,0)\to\bG_m^4$ for different tangent vectors are exactly the
surfaces described above.

\subsubsection*{Case F}
\begin{figure}[htbp!]
  \centering
  \includegraphics[width=274pt]{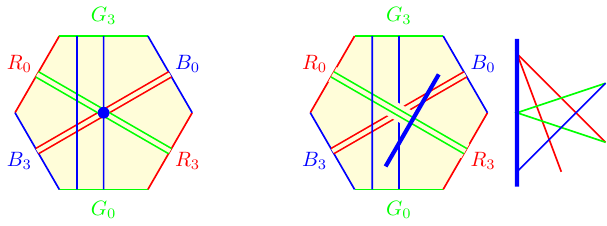} 
  \caption{Degeneration for case F}
  \label{fig-F}
\end{figure}

The configuration is shown in Fig.~\ref{fig-F}. There are five curves
passing through the same point.
It is not log canonical if $c>\frac25$. 
A one-parameter degeneration is resolved by a single blowup of a
point in the central fiber, which becomes $Y_0=\Bl_p\Sigma \cup
\bP^2$. A similar computation as above shows that the moduli space of
the lines in the right panel of Fig.~\ref{fig-F}, not all of
which pass through the same point, is
$(\bA^5/\bA^2)/\bG_m\simeq\bP^2$. 

\subsubsection*{Case G}
Two curves in the same pencil coincide, e.g. $R_1=R_2$. The pair
$(\Sigma, \frac12 D)$ is still log canonical.

\subsubsection*{Case H}
Three curves from the three pencils, e.g. $R_1,G_1,B_1$ pass through
the same point. The pair $(\Sigma, \frac12 D)$ is still log canonical,
and even stronger: it is klt.

\medskip
We also include the following codimension~$2$ degeneration in which a
new type of irreducible components of $Y$ appears.

\subsubsection*{Mixed case EF}

In  case $E$, suppose that $5$ of the $6$ lines in the $\bP^2$ pass
through the same point $q$. Again, this is not a log canonical
configuration if $c>\frac25$, and a generic one-parameter degeneration
with this limit is resolved by one additional blowup of the $3$-fold
at $q$.
For $\frac13<c<\frac12$ this produces the central fiber
$Y_0 = \Bl_p\Sigma\cup \bF_1\cup \bP^2$. For $c=\frac12$ it becomes
$Y_0 = (\bP^1\times\bP^1)\cup \bF_1\cup \bP^2$. 

\begin{figure}[htp!]
  \centering
  \includegraphics[width=295pt]{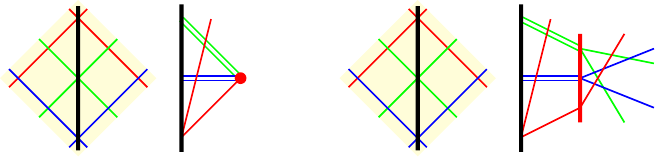} 
  \caption{Degeneration for case EF}
  \label{fig-EF}
\end{figure}

\begin{remark}
  One may ask why we do not consider the configurations in the cases E
  and EF when $6$ lines in the $\bP^2$ on the right pass through a
  common point. The answer is that they do not appear in the
  degenerations $Y\xrightarrow{f} (\Delta,0)\xrightarrow{g}
  \oM(\frac12)$ where the classifying map $g$ is unramified over
  $0\in\Delta$.

  They \emph{do} appear if $g$ if ramified to order $d$
  over $0$.
  (Note that $g(0)$ is a smooth point of $\oM(\frac12)$).
  In that case one has to blow up the central fiber $d$
  times, after which the intermediate ruled surfaces are contracted
  one by one when passing to the canonical model. The picture is quite
  similar to Miles Reid's ``pagoda'' \cite{reid1983minimal-models}.
  For the same reason, in case F we do not consider the case when $5$
  lines pass through a common point.
\end{remark}

\begin{remark}\label{rem: extra-Burniat}
Degenerations of type E can be deformed by smoothing the three pairs of branch lines on the $\pp^2$ component to conics  tangent to the double locus (represented by the black line). These deformations depend on 6 moduli:  3 parameters for each conic tangent to the black line at an assigned point minus the dimension of the group of automorphisms of $\pp^2$ inducing the identity on the black line, which is 3. Similarly, the pairs of fibers of the same color on the $\pp^1\times \pp^1$ component can be smoothed to curves of type $(1,1)$ tangent to the black line. This gives 3 additional moduli (one for each color), so we have an irreducible  9-dimensional family of deformations.

 By \cite[Ex.~1.17]{alexeev2012non-normal-abelian}, for  each such deformation $Y$  there is a $\Z_2^2$-cover $X\to Y$ with $X$ a stable surface that is a deformation of a limit of smooth Burniat surfaces. By dimension reasons all these surfaces are contained in a component of the moduli space of stable surfaces different from $\oM\bur\slc$. 
 
 The situations is similar  for minimal degenerations of type C and D (cf. Figure \ref{fig-mindegs}): we get three additional parameters for type C and  two for type D.. 

\end{remark}

\bigskip

\subsection{The moduli spaces $\oM(\frac13)=\oM\tor$,  $\oM(\frac25)$ and $\oM(\frac12)$}
\label{sec:25}

\begin{theorem}
  For any $\frac14 < c \le \frac13$, 
  the toric family $f\colon \cY\to\oM\tor$ is
  a family of stable pairs
  $\big(Y,\sum_{i=0,3} \frac12 (R_i+G_i+B_i) + \sum_{i=1,2} c(R_i+G_i+B_i)\big)$.
  The fibers are distinct. $\oM(\frac25)$ is projective and the subset
  for which $Y=\Sigma$ is open and dense.
\end{theorem}
\begin{proof}
  By Lemma~\ref{lem:toric-very-ample} and
  Proposition~\ref{prop:family}, 
  the divisors $D\bry$ and $D\inr$ are
  relative Cartier divisors,
  $K+\frac12 D\bry + cD\inr$ is $\Q$-linearly equivalent to $(-1+\frac12 + 2c)K$, hence it is ample for
  $c>\frac14$. By looking at the degenerations one observes that for
  $c\le \frac13$ the pair is log canonical away from the toric
  boundary, i.e. no nontoric degenerations occur.

  We now check that the fibers in the family $\cY\to\oM\tor$ of
  Section~\ref{sec:toric-family} are pairwise non-isomorphic. If
  $Y=\Sigma$ then this amounts to showing that the pair
  $(Y,R_i,G_i,B_i)$ uniquely defines the $6$-tuple
  $(r_1,r_2,g_1,g_2,b_1,b_2)$ in the Equation~\eqref{eq:f}, up to the
  action of $\bG_m^2$. This $6$-tuple defines a unique embedding
  $Y\subset V_P$ since the projective coordinates of $V_p$ are the
  monomials in $r_i$, $g_i$, $b_i$ with the boundary curves.  But this
  is obvious, since these six coefficients are the coordinates of the
  intersections of the curves $R_i$, $G_i$, $B_i$ with the boundary
  curves.

  For a general $Y=\cup Y_k$ one goes through a similar argument for
  the irreducible components $Y_k$. All the types of irreducible
  components already appear in Figs.~\ref{fig-mindegs} and
  \ref{fig-maxdegs}.
\end{proof}

\begin{definition}
  Motivated by  this theorem, we set $\oM(\frac13):= \oM\tor$. 
\end{definition}

The locus in $\oM(\frac13)$ of the points for which the pairs
$(Y,\frac12 D\bry + \frac12 D\inr)$ are not semi log canonical is six
curves $\cC_k$ corresponding to the nontoric degeneration of type F,
when $5$ of the curves pass through a single point. Each of these
curves is $\cC_k\simeq\bP^1$: we include the cases when the remaining
sixth curve goes to the boundary. These six curves intersect at the
origin of the torus $1\in\bG_m^4\subset\oM(\frac13)$, corresponding to
 case E.

\begin{definition}
  Let $\rho_1\colon \oM(\frac25)\to\oM(\frac13)$ be the blowup at the
  origin $1\in\bG_m^4\subset\oM(\frac13)$ with  exceptional divisor
  $E\simeq\bP^3$. Let $\cY(\frac25)\to\oM(\frac25)$ be the blowup in
  \begin{displaymath}
    \cY \times_{\oM(\frac13)} \oM(\tfrac25)
  \end{displaymath}
  of $Z\simeq E$, the preimage of the point of the point
  $p = R_1\cap R_2\cap G_1\cap G_2\cap B_1\cap B_2$. 
\end{definition}

Thus, the fibers of
  $\cY(\frac25)\to\oM(\frac25)$ are the same as the fibers of
  $\cY(\frac13)\to\oM(\frac13)$ outside of $E$, and they are
  isomorphic to $\Bl_p\Sigma\cup\bP^2$ over $E$, as described in case E.

\begin{theorem}
  For any $\frac13<c\le \frac25$, 
  the family
  $\cY(\frac25) \to\oM(\frac25)$
  is a family of stable pairs
  $\big(Y,\sum_{i=0,3} \frac12 (R_i+G_i+B_i) + \sum_{i=1,2} c(R_i+G_i+B_i)\big)$.
  The fibers are distinct. $\oM(\frac25)$ is projective and the subset
  for which $Y\cong \Sigma$ is open and dense.
\end{theorem}
\begin{proof}
  This follows immediately from the computations in the case E
  above. The variety $\oM(\frac 25)$ is a blowup of the projective
  variety $\oM(\frac13)$, so it is projective. According to the
  computation in the case E, the additional fibers are parameterized
  by $E=\bP^3$, so all the fibers in this family are distinct.
\end{proof}

On $\oM(\frac25)$ the preimages of the six curves are disjoint and
lie in the smooth locus. Let us denote them by $\cC_k$ again.

\begin{definition}
  Let $\rho_2\colon\oM(\frac12)\to\oM(\frac25)$ be the blowup along
  $\cup \cC_k$ with  exceptional divisors $E_k$, each of them a
  $\bP^2$-bundle over~$\cC_k$. Let $\cY(\frac12^-) \to\oM(\frac12)$
  be the blowup in
  \begin{displaymath}
    \cY(\tfrac25) \times_{\oM(\frac25)} \oM(\tfrac12)
  \end{displaymath}
  of the preimages of the points of intersection of $5$ lines,
  isomorphic to $\cup E_k$.
\end{definition}

  The fibers are described in the cases F and EF above:
  The fibers of
  $\cY(\frac12^-)\to\oM(\frac12)$ are the same as the fibers of
  $\cY(\frac25)\to\oM(\frac25)$ outside of $\cup E_k$.
  They are isomorphic to $\Bl_p\Sigma\cup\bP^2$ over
  $\cup E_k$ outside of the strict preimage of $E$. 
  And the fibers are isomorphic to $\Bl_p\Sigma\cup\Bl_q\bP^2\cup\bP^2$ over
  $\cup E_k$ intersected with the strict preimage of $E$. 
  Further, let $\cY(\frac12^-)\to\cY(\frac12)$ be the contraction
  defined by $K_{\cY(\frac12^-)} + \frac12\cD\bry + \frac12\cD\inr$,
  relative over $\oM(\frac12)$. This replaces $\Bl_p\Sigma$ by
  $\bP^1\times\bP^1$ as described in the cases E and EF.

  \begin{theorem}\label{thm:Y-family}
  For any $\frac25<c< \frac12$ (resp. for $c=\frac12$)
  the family
  $\cY(\frac12^-) \to\oM(\frac12)$ (resp. $\cY(\frac12) \to\oM(\frac12)$)
  is a family of stable pairs \newline
  $\big(Y,\sum_{i=0,3} \frac12 (R_i+G_i+B_i) + \sum_{i=1,2} c(R_i+G_i+B_i)\big)$.
  The fibers are distinct. \linebreak
  $\oM(\frac12)$ is projective and the subset
  for which $Y=\Sigma$ is open and dense.
\end{theorem}
\begin{proof}
  The proof is immediate from the description in the cases E, F,
  EF. The contraction $\cY(\frac12^-)\to\cY(\frac12)$ exists by
  applying the Basepoint-Free Theorem
  \cite[Thm.~3.24]{kollar1998birational-geometry} over the base
  $\oM(\frac12)$.
\end{proof}

\subsection{All degenerations of pairs, summarized}
\label{sec:irred-comps} 

As a summary, in Table~\ref{tab:Y-irred-comps} we list the types or
irreducible components of the stable pairs
$(Y, \sum_{i=0}^3\frac12(R_i+G_i+B_i))$ that appear in the family over
$\oM(\frac12)$.
For completeness, we also include the minimal degenerations B, G, H in
which $Y=\Sigma$ but the $\bZ_2^2$-cover $X$ is not smooth.

\begin{table}[htbp!]
  \centering
  \begin{tabular}[ht]{lllcl}
    \# &Cases & Surface & Vol & Log canonical conditions\\
    \hline
    0&& $\Sigma$ &6\\
    0&B& $\Sigma$ &6& {\color{gray} $(r_0+r_1\le1,\ g_3+r_1\le 1)$}\\
    0&G& $\Sigma$ &6& {\color{gray} $(r_1+r_2\le1)$}\\
    0&H& $\Sigma$ &6& {\color{gray} $(r_1+g_1+b_1\le 2)$}\\
    1&A& $\bP^1\times\bP^1$ &2&($r_0+r_1+r_2\le1$,\ 
                              $g_3+r_1+r_2\le1$)\\
    2&C&$\bP^1\times\bP^1$ &2&($r_3+r_2+b_1\le1$,\
                             $g_0+g_1+r_2\le1$,\\
       &&&&${\color{gray}b_0+b_1\le1,\ b_3+g_1\le1}$)\\                             
    3&D&$\bF_1$ &3& ($r_0+r_1+b_2\le1$,\
                  ${\color{gray}g_3+r_1\le1,\ b_3+b_2\le1}$)\\
    4&CD, D$^2$&$\bP^2$ &1&($r_0+r_1+b_2\le1,\
                       b_0+b_1+g_2\le1$,\\
       &&&&${\color{gray} b_3+b_2\le1,\ r_3+b_1\le1}$)\\
    5&E, EF&$\bP^1\times\bP^1$ &2&($r_1+r_2+g_1+g_2+b_1+b_2\le2$)\\
    6&E&$\bP^2$ &4& ($r_0+g_0+b_0\le1$)\\
    7&F&$\Bl_1\Sigma$ &5&($r_1+r_2+g_1+g_2+b_1\le2$)\\
    8&F, EF&$\bP^2$ &1&($r_0+g_0+b_0+b_2\le1$)\\
    9&EF& $\bF_1$ &3&($r_0+g_0+b_0\le1$,\ $r_1+r_2+g_1+g_2+b_1\le2$,\\
       &&&&${\color{gray}r_1+r_2\le1,\ g_1+g_2\le1}$)\\
    \hline
  \end{tabular}
  \smallskip
  \caption{Irreducible components of stable pairs in $\oM(\frac12)$}
  \label{tab:Y-irred-comps}
\end{table}

The volume of a component $Y_k$ is
\begin{displaymath}
  4 \big( (K_Y+\tfrac12 D\bry+\tfrac12 D\inr)|_{Y_k}\big)^2 = 
  \big( K_X|_{X_k}\big)^2,
  \quad\text{where } X_k = X\times_Y Y_k.
\end{displaymath}
Note that in a degeneration $X=\cup X_k$ the volumes add up to $K_X^2=6$. 

The last column lists the conditions for the pair
$\big(\bP^2,\sum_{i=0}^2 (r_iR_i+g_iG_i+b_iB_i)\big)$, resp.
$\big(\Sigma,\sum_{i=0}^3 (r_iR_i+g_iG_i+b_iB_i)\big)$, from which
this component of $Y$ originates
(as can be seen by the methods of Section~\ref{sec:wha}) to be log canonical, using the
notation of Lemma~\ref{lem:lc-conditions}.
The main inequalities, in
black, are the ones that fail when all $r_i=g_i=b_i=\frac12$,
i.e. lead to non log canonical singularities, and the ones in gray
lead to log canonical singularities. It is clear that the former
inequalities correspond to the double curves of~$Y$. By
observation, the latter ones correspond to the log centers
of the pair which are contained in $D\inr$.

\begin{lemma}\label{lem:2K_Y+D}
  Up to an isomorphism, there are $10$ surfaces $Y$ underlying the
  stable pairs in $\oM(\frac12)$, with the irreducible components:
  $\Sigma$, $3(\bP^1\times\bP^1)$, $2\bF_1$, two types of
  $2(\bP^1\times\bP^1)\cup2\bP^2$, $(\bP^1\times\bP^1)\cup4\bP^2$,
  $6\bP^2$, $(\bP^1\times\bP^1)\cup\bP^2$, $(\Bl_1\Sigma)\cup\bP^2$
  and $(\bP^1\times\bP^1)\cup\bF_1\cup\bP^2$.  For each of them, the
  divisors $K_Y$ and $D=D\tot=2D\hur$ are Cartier, $D$ is ample, and
  $2K_Y+D$ is very ample.
\end{lemma}
\begin{proof}
  In the toric cases these surfaces appear in the minimal cases A, C,
  D and the maximal cases of Fig.~\ref{fig-maxdegs}. For the nontoric
  cases there are additionally cases E, F and EF. It is easy to see
  that the gluings of the irreducible components $Y=\cup Y_k$ are
  unique. 

  For the toric fibers, we checked that 
  $D$ is ample and $2K_Y+D$ is very ample in
  Lemma~\ref{lem:toric-very-ample} and Proposition~\ref{prop:family}:
  both are sections of $\cO(1)$. Thus, it is enough to look at the
  surfaces appearing in the nontoric cases E and EF.

  In case E, $2K_Y+D$ restricted to $\bP^1\times\bP^1$ and to
  $\bP^2$ is $\cO(1,1)$ and $\cO(2)$ respectively. Consider
  $\bP^3\cup\bP^5\subset \bP^6$ with $\bP^3\cap\bP^5=\bP^2$.  The join
  of the Segre embedding $\bP^1\times\bP^1\subset\bP^3$ and the
  Veronese embedding $v_2\colon\bP^2\to\bP^5$ intersecting along a
  conic in $\bP^2$ is our surface $Y$ embedded by $|2K_Y+D|$. Case
  EF is obtained by a further toric degeneration of
  $v_2(\bP^2)\subset\bP^5$ into a union $\bF_1\cup\bP^2$ with
  $(2K_Y+D)|_{\bF_1}=\cO(s_1+2f)$ and
  $(2K_Y+D)|_{\bP^2}=\cO(1)$. Thus, $|2K_Y+D|$ embeds $\bF_1\cup\bP^2$
  into $\bP^4\cup\bP^2\subset\bP^5$ intersecting along a line
  $\bP^4\cap\bP^2=\bP^1$.

\end{proof}

\subsection{Proof of Theorem~\ref{thm-intro:burniat}}
\label{sec:proof-thm2}

We start with an instructive example.

\begin{example}\label{ex:Burniat-limit-A}
  Consider the degeneration described in case A shown in
  Fig.~\ref{fig-A}. Let $\cY\ini=\Sigma\times\Delta\to(\Delta,0)$ be
  the initial family, and let $\sigma\colon\cY\to\cY\ini$ be the
  composition of two blowups with  exceptional divisors $V_1$,
  $V_2$, so that the central fiber of the family $\cY\to(\Delta,0)$
  is $$Y=V_0\cup V_1\cup V_2=\Sigma\cup\Bl_1\bF_1\cup \bF_0.$$

  On $\cY\ini$ the divisors $\cR\ini$, $\cG\ini$, $\cB\ini$ 
  are
  pullbacks from $\Sigma$, so there exist the sheaves $L_\chi$ pulled
  back from $\Sigma$ which together with the branch divisors $\cR\ini$,
  $\cG\ini$, $\cB\ini$ provide the building data for a $\bZ_2^2$-cover
  $\cX\ini\to\cY\ini$, and $\cX\ini$ is a normal variety.
 
  Let $\cX = \cX\ini \times_{\cY\ini} \cY$ and let $\wcX$ be its
  normalization. Further, let $(\Delta',0)\to(\Delta,0)$ be the
  base change of degree~$2$ ramified at~$0$,
  \begin{displaymath}
    \cY' = \cY \times_{\Delta} \Delta', \quad
    \cX' = \wcX \times_{\Delta} \Delta',
  \end{displaymath}
  and let $\wcX'$ be the normalization of $\cX'$. We claim that
  \begin{enumerate}
  \item   $\cY'$ is a normal variety with  central fiber isomorphic
    to $Y=V_0\cup V_1\cup V_2$,
  \item   $\pi'\colon \wcX'\to\cY'$ is a $\bZ_2^2$-cover whose  branch
    divisors  do not contain any components of the fibers of
    $\cY'\to\Delta'$, and 
  \item the restrictions of $\pi'$ to $V_0$, $V_1$, $V_2$ are
    $\bZ_2^2$-covers such that the curves $V_i\cap V_j$ appear in $R$,
    $G$, $B$ exactly as shown in Fig.~\ref{fig-A}.
  \end{enumerate}

  Indeed, the branch divisors for $\cX\to\cY$ are
  $\sigma^*(\cR\ini)$, $\sigma^*(\cG\ini)$, $\sigma^*(\cB\ini)$, and it is
  clear that $\sigma^*(\cR\ini)$ contains $3V_1+2V_2$ and
  $\sigma^*(\cG\ini)$ contains $V_2$. By
  Lemma~\ref{lem:cover-normalization}, the branch divisors for
  the normalized cover $\wcX\to\cY$ contain $V_1$ and $V_2$ instead,
  so that the total branch divisor is reduced.

  Now make the degree~$2$ base change. Part (1) is obvious.
  The branch divisors for $\cX'\to\cY'$ now contain $2V_1$ and
  $2V_2$ respectively. By Lemma~\ref{lem:cover-normalization} again,
  the branch divisors for $\wcX'\to\cY'$ do not contain either,
  proving part (2). Part~(3) follows by Lemma~\ref{lem:double-curve-colors}.
 
  \medskip
  
  The variety $\wcX'$ is normal and the family
  $\wcX'\to\bA^1_s$ is flat. By Lemma~\ref{lem:slc-cover}
  $K_{\wcX'} = \pi^*(K_{\cY'} + \frac12 D\tot)$ is relatively big
  and nef. Moreover, by the Hurwitz formula
  \begin{displaymath}
    K_{\wcX'} + \wcX'_0 =
    \pi^* \left( K_{\cY'} + \tfrac12 D + \cY'_0 \right).
  \end{displaymath}
  Since the central fiber $(\cY'_0, \frac12 D)$ is slc--in fact it has
  simple normal crossings--the pair
  $\left(\cY', \tfrac12 D + \cY'_0 \right)$ is slc by Inversion of
  Adjunction.  By Lemma~\ref{lem:slc-cover} the pair
  $( \wcX', \wcX'_0 )$ is slc. By Adjunction, implies that $\wcX'$ is
  slc and that $\wcX'$ has canonical singularities along $\wcX'_0$.
  Then by the standard construction, already contained in
  \cite{kollar1988threefolds-and-deformations}, its relative canonical
  model over $\Delta'$, obtained by contraction by the linear system
  $|NK_{\wcX'}|$ for some $N\gg0$, provides a family of stable
  pairs. So we have described the stable limit of Burniat surfaces in
  this example.
\end{example}

\begin{lemma}[Colors for the double crossing locus add up]
  \label{lem:double-curve-colors}
  Let $\cY\to(\Delta,0)$ be a flat family over a smooth curve with 
  central fiber $V_1\cup V_2$ such that $\cY$ is normal and it is
  generically smooth along $V_1\cap V_2$.
  
  Let $\pi\colon\cX\to\cY$ be a $\bZ_2^k$-cover with normal $\cX$ such
  that $V_i$ appear in branch divisors $D_{g_i}$,
  $i=1,2$. (Here, $g_i=0$ means that $V_i$ is not in $D\tot$.)
  Let $(\Delta',0)\to(\Delta,0)$ be a degree~$2$ cover ramified
  over~$0$, $\cY' = \cY\times_{\Delta}\Delta'$,
  $\cX' = \cX\times_{\Delta}\Delta'$, and let $\wcX'$ be the
  normalization of $\cX'$.

  Then the total branch divisor of the $\bZ_2^k$-cover
  $\pi'\colon \wcX'\to\cY'$ does not contain the $V_i$, and $V_1\cap V_2$  is
  contained in the branch divisor $D_{g_1+g_2}$ of $\pi'$.
\end{lemma}
\begin{proof}
  For the cover $\cX'\to\cY'$ either $g_i=0$ and  $V_i$ is not in the branch locus  or $g_i\ne 0$ and $V_i$  is contained in $D_{g_i}$  with multiplicity 2. 
 So by Lemma~\ref{lem:cover-normalization}, for the
  normalized cover $\wcX'\to\cY'$ the branch  divisors 
  contain neither  $V_1$ nor $V_2$.

  Generically along $V_1\cap V_2$, the variety $\cY'$ has an
  $A_1$-singularity, and the components $V_1$, $V_2$ of the central
  fiber of $\cY'\to\Delta'$ are not Cartier along $V_1\cap V_2$.
  Let $U$ be a sufficiently small  open neighborhood of
  $V_1\cap V_2$ in $\cY'$ and $\wU\to U$ be the resolution, with 
  exceptional divisor $E\subset \wU$. 

  For the cover $\cZ'=\cX'\times_{\cY'} \wU \to \wU$ the branch
  divisors are the pullbacks from~$U$. So, $E$ appears in $D_{g_1}$
  and $D_{g_2}$ with coefficient~$1$ if $g_1\ne g_2$ or with
  coefficient~$2$ if $g_1=g_2$.  By
  Lemma~\ref{lem:cover-normalization} for the normalized cover
  $\wcZ'\to \wU$, $E$ appears in the branch divisor
  $D_{g_1+g_2}$. So in the restriction of $\wcZ'\to \wU$ to $V_i$ the
  divisor $V_i\cap E$ appears in $D_{g_1+g_2}$. The double locus
  $V_i\cap E$ equals $V_i\cap V_j\cap U$,
  with the same neighborhoods in $V_i$,
  and the normalized
  covers $\wcX'$, $\wcZ'$ agree.  This proves the statement.
\end{proof}

We are now ready to prove Theorem~\ref{thm-intro:burniat}. We will
give two proofs: (1) by analyzing the stable limits of one-parameter
families, and (2) by constructing families of stable Burniat
surfaces over an open cover of the moduli space $\oM(\frac12)$. 

\begin{proof}[First proof of Theorem~\ref{thm-intro:burniat}]
  Consider a one-parameter family $\cX\to \Delta\setminus 0$ of
  Burniat surfaces. After a finite base change, it can be realized as
  a $\bZ_2^2$-cover $\cX\to\cY$ of a family of stable pairs
  $(\cY, \frac12\cR+\frac12 \cG + \frac12 \cB)$ over a punctured
  one-dimensional base $\Delta\setminus 0$. To simplify the notation,
  let us assume that we have made this reduction. 

  Since $\oM(\frac12)$ is complete, 
  we have an extension $\Delta\to\oM(\frac12)$
  and the pullback of the family of Theorem~\ref{thm:Y-family} over
  $\oM(\frac12)$ gives a family over $\Delta$ whose
  central fiber $(Y,\frac12 R+\frac12 G+\frac12 B)$ is
  stable. 
  We claim that the stable
  limit of the Burniat family is the $\bZ_2^2$-cover of $Y$.
  In particular, the stable limit of Burniat surfaces is uniquely defined
  by the image of $0\in\Delta$ in $\oM(\frac12)$. Since $\oM(\frac12)$
  is normal and the moduli space of smooth Burniat surfaces is
  $M^0/ \Gamma$, this implies that the normalization of the
  closure of the moduli space of Burniat surfaces in the moduli
  space of stable surfaces is $\oM(\frac12)/ \Gamma$, thus proving
  Theorem~\ref{thm-intro:burniat}.

  The proof of the above statement is the same as in
  Example~\ref{ex:Burniat-limit-A}, with minor changes. We start with
  the normal variety $\cY\to\Delta$. The solution to the fundamental
  relations exists over $\Delta\setminus 0$. Then,
  denoting by $\cR,\cG,\cB$ the closures of the divisors over
  $\Delta\setminus 0$, 
  there exist
  divisorial sheaves $L_1,L_2,L_3$ on $\cY$ such that in the class
  group $\Cl(\cY)$ one has
  \begin{displaymath}
    2L_1=\cG + \cB + \sum V_{1,k}
  \end{displaymath}
  for some components $V_{1,k}$ of the central fiber $Y$, and
  similarly for $L_2$ and $L_3$.  We make a $2:1$ base change
  $(\Delta',0)\to(\Delta,0)$,
  so that the pullbacks of $V_{1,k}$ are $2V_{1,k}$, 
  and find the building data with 
  reduced total branch divisor. As in Example \ref{ex:Burniat-limit-A}, this removes the vertical components from the total  branch divisor. 
  If $\wcX'\to\cY'$ is the
  $\bZ_2^2$-cover for these data then by Lemma~\ref{lem:slc-cover} the
  pair $(\wcX',\wcX'_0)$ has slc singularities and $K_{\wcX'}$ is
  relatively ample. By adjunction, it follows that
  \begin{enumerate}\label{enu:CM}
  \item $\wcX'$ has canonical singularities near $\wcX'_0$,
    thus both are Cohen-Macaulay,
  \item the central fiber $\wcX'_0$ has slc singularities, and 
  \item $K_{\wcX'}$ is relatively ample.
  \end{enumerate}
  Therefore, $\wcX'_0$ is the stable limit of
  $\cX\to\Delta\setminus0$. The induced map $\wcX'_0\to\cY$ is a
  $\bZ_2^2$-cover between deminormal varieties with slc singularities.
  The restrictions of this map to the irreducible components $Y_i$ of
  $Y$ are also $\bZ_2^2$-covers. In all cases, excluding case 4 of
  Table~\ref{tab:Y-irred-comps}, ``colors'' of curves in the double
  crossing locus are uniquely determined by the fact that the
  fundamental relations~\eqref{eq:fundrel} on the irreducible
  components $Y_i$ have a solution. When some component $Y_i$ is of
  type 4, we look at a two-step degeneration, from a simpler limits
  where the components are $\bP^1\times\bP^1$ or $\bF_1$. This implies
  that the colors are uniquely determined as well.

  So the stable limit of Burniat surfaces is uniquely determined by
  the stable limit $(Y,\frac12 R+\frac12 G+\frac12 B)$, as claimed.
  The normalization map is a bijection by Lemma~\ref{lem:thm2-stronger}.
\end{proof}

For another proof, we construct a family of stable Burniat
surfaces over open sets covering $\oM(\frac12)$, after appropriate
$2^n:1$ base changes. 

\begin{proof}[Second proof of Theorem~\ref{thm-intro:burniat}]
Consider the family $f\colon\cY(\frac12)\to\oM(\frac12)$ of
Theorem~\ref{thm:Y-family}. Over the open subset
\begin{displaymath}
  U = \oM(\tfrac12) \setminus \left[ \text{(toric boundary)}
 \cup( \cup_{k=1}^6 E_k)\right] =
  \bG_m^4\setminus \cup_{k=1}^6\cC_k  
\end{displaymath}
the morphism $f$ is smooth,
every fiber is isomorphic to $\Sigma$ and we have a family of
effective divisors $R_i,G_i,B_i$, $i=0,1,2,3$.  The variety
$\cY(\frac12)$ is normal and the complement $\oM(\frac12)\setminus U$
is a union of divisors. 

Define the Weyl divisors $\cR_i,\cG_i,\cB_i$ on $\cY(\frac12)$ to be
the closures of the divisors over $f\inv(U)$.
Consider the
divisors $\cR~+\cG$, $\cG+\cB$, $\cB+\cR$. (As before, $\cR=\sum\cR_i$,
$\cG=\sum\cG_i$, $\cB=\sum\cB_i$.)  Their restrictions to every fiber
over a point in $U$ are $2$-divisible. Since $\Cl(U)=0$, 
there exist divisorial sheaves $L_1,L_2,L_3$ on $\oM(\frac12)$ with
\begin{displaymath}
  2L_1 = \cG+\cB + \sum\cV_{1,k}
  \text{ in } \Cl\big( \cY(\tfrac12) \big)
\end{displaymath}
for some Weil divisors $\cV_{1,k}$ supported on the boundary
$f\inv\big(\oM(\frac12)\setminus U\big)$, and similarly for $L_2$ and $L_3$.

Let $s\in \oM(\frac12)\setminus U$ be a point on the boundary. Suppose
that it is an intersection of the (finitely many) boundary divisors
$H_n$. Pick a sufficiently small affine neighborhood $s\in
W\subset\oM(\frac12)$. There exists a finite cover $\mu\colon V\to W$ with
normal $V$ such that
\begin{enumerate}
\item $\mu$ is branched along each $H_n$ with multiplicity divisible by $2$,
\item the base changed family
  $\cY_{V} = \cY(\frac12)\times_{\oM(\frac12)} V \to V$
  is normal.
\end{enumerate}
Indeed, outside of the $8$ singular torus-fixed points (see
Corollary~\ref{cor:M-sings}) the boundary of $U$ is normal-crossing, and
near the toric boundary we can take the toric base change
corresponding to the embedding $2N\to N$ of the cocharacter lattices.

Since $\Cl(U)=0$, it follows that $\Cl(\oM(\frac12))$ is generated by
$H_n$. On the cover the multiplicities of $H_n$ double. So 
the pullbacks of $\cR+\cG$, $\cG+\cB$, $\cB+\cR$ 
to $\cY_{V}$ are divisible by~$2$ 
and thus satisfy the fundamental
relations \eqref{eq:fundrel} for a $\bZ_2^2$-cover
$\pi_V\colon\cX_V\to\cY_V$.
The restriction of this family to a generic one-parameter subfamily is
the same as described in the first proof. So $\pi_V$ is a family of
stable surfaces  and it induces a classifying morphism
$V\to \oM\bur\slc$. Clearly, it descends to $W\subset\oM(\frac12)$
and for different open neighborhoods $W_1,W_2$ the map is the same on
$W_1\cap W_2$. So we get a morphism $\oM(\frac12)\to\oM\slc\bur$.

Finally, we have to divide by the relabeling group for different choices
of labeling for smooth Burniat surfaces, giving a well-defined
morphism $\oM(\frac12)/\Gamma \to \oM\bur\slc$.  We note
here that the natural $\Gamma$-action on the fan $\fF$ gives
the action on $\oM(\frac13)$, and that the unions of centers of the
blowups $\rho_1$, $\rho_2$ in Section~\ref{sec:25} are invariant, giving a
natural $\Gamma$-action on $\oM(\frac12)$.

The
morphism $\oM(\frac12)/\Gamma \to \oM\bur\slc$
is finite, it is a bijection on a dense open subset, and the
source is a normal variety. It follows that the normalization of an
irreducible component of $\oM\bur\slc$ is indeed
$\oM(\frac12)/\Gamma$.
\end{proof}

\begin{lemma}\label{lem:2K_X}
  Let $X$ be a degenerate Burniat surface and let $\pi\colon X\to Y$
  be the corresponding $\bZ_2^2$-cover.  Then, denoting
  $D=D\tot=2D\hur$, the linear system $|2K_X|$ coincides with
  $\pi^* |2K_Y+D|$, it is base point free and maps $X$ to $Y\subset\bP^6$.
\end{lemma}
\begin{proof}
  The surfaces $X$ and $Y$ are slc. By Lemma~\ref{lem:2K_Y+D}, $Y$ is
  Gorenstein and $D$ is Cartier.  By Lemma~\ref{lem:slc-cover} we have
  the equality of line bundles $2K_X = \pi^*(2K_Y+D)$.

  For $i=1,2$ one has $h^i(2K_Y+D)=h^{2-i}(-K_Y-D)$, and the latter is
  zero by \cite[Prop.~3.1]{liu-rollenske2014} since $K_Y+D$ is ample
  by Lemma~\ref{lem:2K_Y+D}.  Similarly, $h^i(2K_X)=0$ for
  $i=1,2$. Thus, $h^0(2K_X)=\chi(2K_X)$ and
  $h^0(2K_Y+D)=\chi(2K_Y+D)$.

  Both $X$ and $(Y,\frac12 D)$ are flat limits of 
  surfaces and the Euler characteristic is locally constant in flat
  families. Thus, $\chi(2K_X)=\chi(2K_Y+D)=7$ since this holds
  generically. Thus, $h^0(2K_X)=h^0(2K_Y+D)$, which implies $|2K_X| =
  \pi^*|2K_Y+D|$.
  The last part follows by Lemma~\ref{lem:2K_Y+D}.
\end{proof} 

The following Lemma implies the last claim of Theorem~~\ref{thm-intro:burniat}.

\begin{lemma}\label{lem:thm2-stronger}
  The degenerate Burniat surfaces corresponding to different points of
  $\oM(\frac12) / \Gamma$ are non-isomorphic. In other
  words, the map from $\oM(\frac12) / \Gamma$ to the 
  moduli space of stable surfaces is a bijection to the closure of $M\bur$.
\end{lemma}
\begin{proof}
  By Theorem~\ref{thm:Y-family}, the fibers
  $(Y,\sum_{i=0}^3\frac12(R_i+G_i+B_i))$ over different points of
  $\oM(\frac12)$ are non-isomorphic.  By the previous
  Lemma~\ref{lem:2K_X}, if $-K_Y$ is very ample then we one can
  recover the pair $(Y,\sum_{i=0}^3\frac12(R_i+G_i+B_i))$ from $X$,
  since the $\bZ_2^2$-cover is intrinsic: it is the bicanonical
  map. By Lemma~\ref{lem:2K_Y+D} this is the case in all cases except
  when $Y=(\bP^1\times\bP^1)\cup\bF_1\cup\bP^2$. In the latter case,
  $|2K_X|$ and $|-K_Y|$ map $X$ and $Y$ to
  $(\bP^1\times\bP^1)\cup\bP^2$. But we observe that the middle
  component $\bF_1$ in the pair $(Y,\sum_{i=0}^3\frac12(R_i+G_i+B_i))$
  is unique and does not vary in moduli.
\end{proof}

\subsection{Degenerate Burniat surfaces} 
\label{sec:burniat-fibers}

\begin{theorem}
  The boundary of the moduli space $M\bur$ of smooth Burniat surfaces
  in the compactification $\oM\bur\slc$ consists of $8$ divisors
  corresponding to the degenerations A, B, C, D, E, F, G, H of
  Sections~\ref{sec:mindegs} and \ref{sec:nontoric}.
\end{theorem}
\begin{proof}
  These are all the minimal degenerations modulo $\Gamma$
  and for each of them we found an irreducible $3$-dimensional family
  of pairs in $\oM(\frac12)$ and of their $\bZ_2^2$-covers. Thus, the
  closure of each set is a divisor in
  $\oM\bur\slc = \oM(\frac12)/ \Gamma$. Outside of the
  union of these divisors one has $Y=\Sigma$ and the curves $R_i$,
  $G_i$, $B_i$ are in general position. So the covers are smooth
  Burniat surfaces. 
\end{proof}

We now describe the degenerate Burniat surfaces for a general point in
each of these divisors. 
\begin{remark}
All the surfaces in $\oM\bur\slc$ satisfy $h^1(\OO)=h^2(\OO)=0$, since they are stable limits of smooth surfaces and slc singularities are Du Bois (\cite{kollar2010-dubois}).  
It is also possible to double check this vanishing directly in all cases, as is done  in Examples 3.5 and 3.6  of \cite{alexeev2012non-normal-abelian} for the last degeneration to the right in Figure \ref{fig-maxdegs} and  for the  degeneration of type E.
\end{remark}

\begin{warning}
  It frequently happens that the components $X_k$ of a stable Burniat surface $X$ are not $S_2$, even though the
  surface $X=\cup X_k$ is. This happens when for some point $p\in X_k$
  the stabilizer group for the $G$-action in $X$ is bigger than 
  the stabilizer group for the corresponding $S_2$ surface
  $X_k^\sigma$. The stabilizer group $G_p$ is generated by the
  elements $R,G,B\in\bZ_2^2$ for the preimages of the curves $R_i$, $G_i$, $B_i$ (in
  $X$, resp. in $X_k$ only) containing $p$. When the two stabilizers
  are different, $X_k$ is obtained from $X_k^\sigma$ by gluing several
  points. 
\end{warning}

For each irreducible component $Y_k$ of a fiber $Y$ of the universal
family $\mathcal Y\to \oM(\tfrac12)$, the $S_2$-fication
$X_k^{\sigma}$ of the $\bZ_2^2$-cover $X_k\to Y_k$ is given by the
procedure described in \cite{alexeev2012non-normal-abelian}, which we
reviewed in Section~\ref{subsec:covers}.  The geometric characters of
$X_k^{\sigma}$ can be all recovered from the branch data: the
canonical class, and therefore $K^2$ and the Kodaira dimension, can be
computed using Lemma \ref{lem:slc-cover}, and the cohomology of the
structure sheaf can be computed using the decomposition
$\OO_Y\bigoplus(\oplus_{i=1}^3L_i\inv)$ of its direct image. (Note
that since all the $Y_k$ are simply connected and therefore, a fortiori, have torsion free Picard group,  their decomposition is
uniquely determined by the branch data). Finally the slc singularities
that can occur in our situation have been analyzed in
\cite{alexeev2012non-normal-abelian}.

\subsubsection*{Case A} (Fig.~\ref{fig-A}).  In the general case, namely when all the lines are
distinct, each component is a smooth bielliptic surface (so
$K^2=p_g=0$, $q=1$) and the Albanese pencil is the pull back of the
ruling of $\pp^1\times\pp^1$ that contains 2 pairs of  branch lines
in different branch divisors  $R$, $G$ or $B$.
Two components are
glued transversally along a smooth elliptic curve. All three
components meet at one point, which is a degenerate cusp of $X$.
 
Another description, useful in understanding the degenerations, is as
follows.  For the general case, consider three elliptic curves $E_1$,
$E_2$ and $E_3$, and on each $E_k$ a translation $\tau_k$ by a point
of order 2 and a rational involution $\si_k$. Let $\si_k'$ be the
involution induced by $\si_k$ on $E'_k:=E_k/\tau_k$. Take
$X_k:=(E_{k+1}\times E'_{k+2})/\bZ_2$, where $\Z_2$ acts on $E_{k+1}$
via $\tau_{k+1}$ and on $E'_{k+2}$ via $\si'_{k+2}$ (the index $k$
varies in $\bZ_3$). The surfaces $X_k$ and $X_{k+1}$ are glued along a
curve isomorphic to $E'_{k+2}$, which on $X_k$ is a fiber of the
Albanese pencil $X_k\to E'_{k+1}$ and on $X_{k+1}$ is half of a fiber
of the rational pencil $X_{k+1}\to E'_k/\si'_k=\pp^1$.
 
Letting two lines in the same branch divisor coincide
corresponds to degenerating one of the $E_k$ to a cycle of two
rational curves. Letting two lines that are in different branch
divisors on one component coincide corresponds to degenerating one of
the $E_k$ to a nodal rational curve. At most  three degenerations of this
type can occur at the same time.

This surface appears very nicely as a degeneration of Burniat surface
in the form given by Inoue \cite{inoue1994some-new-surfaces}, with the
parameter $\lambda\to 0$ or $\infty$. 

\subsubsection*{Case B}  The surface $X$ is non-normal, with
singularities of types $2'.1$, $2'.2$, $3'.2$, $3'.4$, and $4''.6$ in
Tables 2 and 3 of \cite{alexeev2012non-normal-abelian}. The
normalization is a (non-minimal) properly elliptic surface with 2 $A_1$ singularities. 

\subsubsection*{Case C} (Fig.~\ref{fig-C}).
In the general case, the  surfaces $X_1^{\sigma}$, $X_2^{\sigma}$ and $X_3^{\sigma}$, are
singular Enriques surfaces.  The surfaces $X_k$ meet transversally
at one point $p_0$ which is smooth for all of them, so $X$ has a
degenerate cusp there.  Two components $X_k$ and $X_{k+1}$ are glued
along a rational curve with a node $p_{k+2}$. At $p_{k+2}$ there is
additional gluing and the surface such that $p_{k+2}$ lies on 3 lines
in the same branch divisor is not $S_2$ there.

\subsubsection*{Case D} (Fig.~\ref{fig-D}).  Each surface $X_k^{\sigma}$ is a singular  properly
elliptic surface with 2$D_4$ and one $\frac 14(1,1)$ singularities and with $h^1(\OO)=h^2(\OO)=0$.  The elliptic fibration
is given by $|2K|$ and it is the pullback of the ruling of
$\bF_1$. The two components are glued along a 
rational curve with two nodes $p_1$, $p_2$, where there is an
additional gluing. Each component is not $S_2$ at one of the points
$p_k$ (the one with three branch lines of the same color going through the point).

\subsubsection*{Case E} (Fig.~\ref{fig-E}). In this case we have $X_i=X_i^{\sigma}$ for $i=1,2$. The surface $X_1$, the
double cover of $Y_1=\pp^1\times \pp^1$ is a singular  del Pezzo
surface with $K^2=2$,
with 6 $A_1$ singularities.  It is the quotient of $\pp^1\times\pp^1$ by the diagonal action of $\bZ_2^2$. The second component $X_2$ is a
degenerate Enriques surface.
 The two surfaces are glued along a curve with 6 nodes that is the union of 4 smooth rational curves.  If we let two of the pairs of lines in the same branch divisor
coincide on $Y_2$, then $X_2$ becomes reducible  and the normalization  is the union of
two quadric cones. 

\subsubsection*{Case F} (Fig.~\ref{fig-F}). Both
components are $S_2$.
The component $X_1$ which is the cover of the
blow up $Y_1$ of $\Si$ at one point has $K^2=2$, $h^1(\OO)=1$,
$h^2(\OO)=0$.  It is not normal along the preimages of the  two double lines in the branch locus, where it has  double crossings points. The normalization of $X_1$ is a ruled surface with $h^1(\OO)=1$, whose  Albanese pencil is induced by the  pencil  that has no fiber contained in the branch locus with multiplicity 2 (the blue one in the picture).  
The second component $X_2$ is the same  degenerate Enriques surface in case E, but in this case $X_1$ and $X_2$ are glued along the union of two
rational curves meeting transversally at two points.

\subsubsection*{Case G}
The surface becomes non-normal, with double crossings singularities. 
The normalization $X^{\nu}$   is a non minimal  bielliptic surface. In fact if, say, $G_1=G_2$, then contracting  $R_0$ and $R_3$ gives  $\pp^1\times\pp^1$ with the same configuration of lines as in case A, so the cover is a bielliptic surface $\bar X$ and $X^{\nu}$ is the blow up of $\bar X$ at two points. 

\subsubsection*{Case H}
The surface acquires a $\frac14(1,1)$ singularity, with  desingularization  a Burniat surface with $K^2=5$.

\medskip

The cases above involve all of the irreducible components $Y_k$ of
Table~\ref{tab:Y-irred-comps} except for \#4 and \#9. We now describe
the covers in these two cases.

\subsubsection*{\#4} $X_k$ is a singular  del Pezzo surface with $K^2=1$ with 2$D_4$ singularities over the point where three lines of the same color meet and a $\frac14(1,1)$ point over the point where three lines of three distinct colors meet. It
is glued to the neighboring components along two rational curves with a
node.

\subsubsection*{\#9} $X_k$ is a non-normal  surface with $K^2=-4$. 
It has two irreducible components that are del Pezzo surfaces with $K^2=6$ and 2$A_1$ singularities.
has double crossing singularities along two disjoint rational curves.

\begin{remark}
  In Lemma~\ref{lem:cones} we listed the cones of the fan $\fF$ modulo
  $\Gamma$, i.e. the toric degenerations of the pairs
  $(Y,\frac12 D)$. There are $29$ of them. Adding non-toric
  degenerations, there are EF, BF, and then all the possible subcases
  of these $31$ cases obtained by adding some combinations of the G and H
  degenerations. It does not seem practical to list all of these
  possibilities here.
\end{remark}

\begin{remark}
  Although the space $\oM\bur$ which we constructed is irreducible, in
  the larger space of stable surfaces there are definitely other
  irreducible components intersecting $\oM\bur$. For example, in the
  degeneration E of Fig.~\ref{fig-E}
  the pairs of lines on $\bP^2$ can be deformed to
  conics tangent to the double locus. Similarly, the
  three divisors of type $(1,1)$ on $\bP^1\times\bP^1$ can be
  smoothed, keeping them tangent to the double locus. Since the
  induced $\bZ_2^2$ covers of the double curve $\bP^1$ have the same
  normalization, the covers can be glued. This gives a family of
  dimension 12. Many of the other degenerations produce other
  irreducible components in the moduli of stable surfaces.
\end{remark}

\section{An alternative description with weighted line arrangements}
\label{sec:wha}

In the Campedelli case, the compactification of the moduli space of
pairs \linebreak
$(\bP^2,\sum_{i=1}^7 \frac12 D_i)$ is a special case of a more general
situation for the weighted hyperplane arrangements
$(\bP^{r-1}, \sum_{i=1}^n d_i D_i)$ for some fixed weights
$0<d_i\le 1$, considered in \cite{alexeev2015moduli-weighted}. The
theory becomes rather trivial in this particular case.

In the Burniat case, the pairs
$\big(\Bl_3\bP^2, \sum_{i=0}^3 \frac12 (R_i+G_i+B_i)\big)$ arise from a
configuration of $9$ lines in $\bP^2$ (in two ways, related by a
Cremona transformation). The construction of the
compactified moduli spaces can be done as an application.

We sketch this alternative way here, omitting some details.

\subsection{Compact moduli space for weighted hyperplane arrangements}
\label{sec:wha-general}

Let $\beta = (b_1,\dotsc,b_n) \in (0,1]^n$ be a fixed
vector and consider the moduli space $M_\beta(r,n)=M_\beta(\bP^{r-1},n)$ of
log canonical pairs $(\bP^{r-1}, \sum b_i B_i)$ where $B_i$ are some
hyperplanes on $\bP^{r-1}$, some of which are allowed to
coincide. This moduli space can be easily constructed as a free
quotient of an open subset of $(\bP^{r-1}{}^\vee)^n$ by a free action
of $\PGL(r)$, since the log canonical pairs as above have trivial
automorphism groups.

\cite{alexeev2015moduli-weighted} constructs a certain projective
scheme $\oM_\beta(r,n)$ together with a family of stable pairs
$(X,\sum b_iB_i)$, called stable weighted hyperplane arrangements
over it, containing $M_\beta(r,n)$ as an open subset. The theory is a
generalization of \cite{hacking2006compactification-moduli} to the
weighted case, in the same way Hassett's moduli spaces $\oM_{g,\beta}$
of weighted stable curves \cite{hassett2003moduli-spaces} generalize
the Deligne-Mumford-Knudsen's spaces~$\oM_{g,n}$. 
The basic tools used in the description are:

(1) The hypersimplex
\begin{displaymath}
  \Delta(r,n) = \{ (x_1,\dotsc,x_n) \mid 0\le x_i\le 1, \ \sum x_i=r \}.  
\end{displaymath}

(2) Matroid polytopes $P\subset\Delta(r,n)$ associated to hyperplane
arrangements \linebreak $(\bP^{r-1}, B_1,\dotsc, B_n)$. They are
defined as follows: for each subset of indices $I\subset \{1,\dotsc,
n\}$ one adds the inequality
\begin{math}
  \sum_{i\in I} x_i \le \codim \cap_{i\in I} B_i.
\end{math}
It is enough to consider flats, i.e. maximal sets $I$ producing the
same linear space $Z=\cap_{i\in I} B_i$. And one can omit the normal
crossing intersections since for them the inequalities follow from
$x_i\le 1$. One way to understand this matroid polytope is that this is
the set of weights $(x_i)\in [0,1]^n$ for which the pair $(\bP^{r-1},
\sum x_iB_i)$ is log canonical.

(3) The weighted hypersimplex, the ``window''
\begin{displaymath}
  \Delta_\beta(r,n) = \{ (x_1,\dotsc,x_n) \mid 0\le x_i\le b_i, \ \sum x_i=r \}.
\end{displaymath}

\smallskip

Then a stable weighted hyperplane arrangement $(X,\sum b_iB_i)$ is
described by a \emph{partial} 
face-fitting cover $\cup P_k$ of $\Delta(r,n)$ by matroid polytopes
$P_k$ intersecting the interior of $\Delta_\beta(r,n)$ and such that
$\cup P_k$ completely covers the window $\Delta_\beta(r,n)$.

Then the irreducible components $X_k$ of $X$ are in a bijection with the
polytopes $P_k\cap\Delta_\beta(r,n)$, and the configurations of the
divisors $B_i$ on $X_k$ can be read off this combinatorial gadget as
well.

\medskip

We now explain how this general theory applies to the Campedelli line
arrangements and sketch an extension of this theory to the Burniat
line arrangements. 

\subsection{Campedelli line arrangements}
\label{sec:wha-campedelli}

In this case the matroid polytope associated to a line arrangement
$(\bP^2, \sum_{i-1}^7 \frac12 D_i)$ is a subset of
$\Delta(3,7) = \{(x_i)\in [0,1]^7 \mid \sum x_i=3\}$ satisfying the
following additional conditions:
\begin{enumerate}
\item If several lines $B_i$, $i\in I$, coincide then $\sum_{i\in
    I}x_i\le 1$. 
\item If several lines $B_i$, $i\in I$, pass through a common point
  then $\sum_{i\in I}x_i\le2$. 
\end{enumerate}
The weighted ``window'' is 
\begin{displaymath}
  \Delta\halfwts(3,7) = 
  \{ (x_1,\dotsc,x_7) \mid 0\le x_i\le \tfrac12, \ \sum x_i=3 \}.
\end{displaymath}

What makes the situation easy is the fact, easily checked, that none
of the equations $\sum_{i\in I}x_i=1$ and $\sum_{i\in I}x_i=2$ cuts
the interior of $\Delta\halfwts(3,7)$. Thus, in
this case each partial matroid tiling consists of a unique matroid
polytope $P$. As in Remark~\ref{rem:campedelli-all-cases}, there are
$36$ such matroid polytopes modulo $S_7$, and $175$ modulo
$\GL(3,\bF_2)$.

The geometric consequence of this combinatorial statement is that for
every stable weighted hyperplane arrangement
$(X,\sum_{i=1}^7\frac12 B_i)$, the underlying variety $X$ is
irreducible and isomorphic to $\bP^2$, which is the same conclusion
that we arrived to in Section~\ref{sec:compact-campedelli} by an
easier method.

\subsection{Matroid tilings for Burniat arrangements}
\label{sec:burniat-tilings}

Consider the initial generic Burniat line arrangement of the $9$ lines
$R_i, G_i, B_i$, $i=0,1,2$, on $\bP^2$ in
Fig.~\ref{fig-burniat-config}. Associated with it is the matroid
polytope in $\Delta(3,9)$ of the $9$-tuples $(r_i,g_i,b_i)$,
$i=0,1,2$, satisfying the inequalities $0\le r_i,g_i,b_i\le1$,
$\sum(r_i+g_i+b_i)=3$, and
\begin{displaymath}
  g_0+g_1+g_2+b_0\le 2, \ b_0+b_1+b_2+r_0\le 2, \ r_0+r_1+r_2+g_0\le2.
\end{displaymath}
We will call it the Burniat polytope and denote by
$\Delta\ubur$. Recall from Lemma~\ref{lem:lc-conditions} that there
are three additional coordinates
\begin{displaymath}
  r_3=g_0+g_1+g_2+b_0-1, \ g_3=b_0+b_1+b_2+r_0-1, \ b_3=r_0+r_1+r_2+g_0-1.
\end{displaymath}
It is easy to see that $\sum_{i=0}^2(r_i+g_i+b_i)=3 \iff
\sum_{i=1}^3(r_i+g_i+b_i)=3$. Indeed, this is equivalent to 
$K+ \sum (r_iR_i+g_iG_i+b_iB_i) \equiv 0$ on
the respective $\bP^2$'s, and this condition is
preserved by the Cremona transformation. 

In terms of these variables the above inequalities become
$r_3,g_3,b_3\le1$. The situation is not totally symmetric, however,
because it is \emph{not} true that $r_3,g_3,b_3\ge0$ on $\Delta\ubur$.
By analogy with the weighted hyperplane case, we define a ``window'',
the weighted Burniat polytope as the subset
\begin{displaymath}
  \Delta\ubur\halfwts =
  \{ 0\le r_i,g_i,b_i \le\tfrac12 \text{ for all } i=0,1,2,3 \}
  \subset \Delta\ubur
\end{displaymath}

We now observe that the inequalities in
Table~\ref{tab:Y-irred-comps} defining the irreducible components of
the degenerate surfaces $Y$ are in fact the inequalities defining
certain matroid subpolytopes $P_k\subset\Delta\ubur$. Moreover, the
degenerations A, B, \dots, H and EF naturally correspond to the
partial matroid tilings of $\Delta\ubur$ -- completely covering the
window $\Delta\halfwts\ubur$ -- which are listed in
Table~\ref{tab:burniat-min-tilings}.
The inequalities in black are those that cut through the interior of
$\Delta\ubur\halfwts$.  
The reason for this will become
clear in the next section.

\begin{table}[htbp!]
  \centering
  \begin{tabular}{lrl}
    Case &Vol& Matroid polytopes\\
    \hline
    A&2& ($r_0+r_1+r_2\le1$,\ $g_3+r_1+r_2\le1$),\\
     &+2& ($g_0+g_1+g_2\le1$,\ $b_3+g_1+g_2\le1$),\\
     &+2& ($b_0+b_1+b_2\le1$,\ $r_3+b_1+b_2\le1$).\\
    B&6& ${\color{gray} (r_0+r_1\le1,\ r_1+g_3\le1).}$ \\
    C&2& ($r_3+r_2+b_1\le1$,\ $g_0+g_1+r_2\le1$,\ 
       ${\color{gray}b_0+b_1\le1,\ b_3+g_1\le1}$),\\  
     &+2& ($g_3+g_2+r_1\le1$,\ $b_0+b_1+g_2\le1$,\ 
       ${\color{gray}r_0+r_1\le1,\ r_3+b_1\le1}$),\\  
     &+2& ($b_3+b_2+g_1\le1$,\ $r_0+r_1+b_2\le1$,\ 
       ${\color{gray}g_0+g_1\le1,\ g_3+r_1\le1}$).\\  
    D&3& ($r_0+r_1+b_2\le1$,\ ${\color{gray}g_3+r_1\le1,\ b_3+b_2\le1}$),\\
     &+3& ($r_3+r_2+b_1\le1$,\ ${\color{gray}g_0+r_2\le1,\ b_0+b_1\le1}$).\\
    E&2& ($r_1+r_2+g_1+g_2+b_1+b_2\le2$), \\
     &+4& ($r_0+g_0+b_0\le1$).\\
    F&5& ($r_1+r_2+g_1+g_2+b_1\le2$),\\
     &+1& ($r_0+g_0+b_0+b_2\le1$).\\
    G&6& {\color{gray} $(r_1+r_2\le1).$}\\
    H&6& {\color{gray} $(r_1+g_1+b_1\le 2).$}\\
    EF&2& ($r_1+r_2+g_1+g_2+b_1+b_2\le2$), \\
         &+3& ($r_0+g_0+b_0\le1$,\ $r_1+r_2+g_1+g_2+b_1\le2,\
           {\color{gray}r_1+r_2\le1,\ g_1+g_2\le1}$),\\
         &+1& ($r_0+g_0+b_0+b_2\le1$).\\
    \hline
  \end{tabular}
  \smallskip
  \caption{Matroid tilings of $\Delta\ubur\halfwts$ for some degenerations}
  \label{tab:burniat-min-tilings}
\end{table}

\subsection{Compactification for the Burniat arrangements} 
\label{sec:burniat-alternative} 

We give a sketch of a second construction of the
compactification, in addition to the one in
Section~\ref{sec:Burniat}. 

As in Section~\ref{sec:wha-general}, let $M\halfwts(3,9)$ denote the moduli
space of log canonical pairs of $\bP^2$ with $9$ lines, which we label
$R_i, G_i, B_i$, $i=0,1,2$. It comes with the compactification
$\oM\halfwts(3,9)$. Let $Z\subset M\halfwts(3,9)$ be the closed subset
of arrangements for which there are three quadruples of lines passing
through a common point, as in Fig.~\ref{fig-burniat-config}:
\begin{displaymath} 
  R_0\cap R_1\cap R_2\cap G_0=p_B, \
  G_0\cap G_1\cap G_2\cap B_0=p_R, \
  B_0\cap B_1\cap B_2\cap R_0=p_G.
\end{displaymath}
Let $\oZ$ be its closure in $\oM\halfwts(3,9)$, with the reduced
scheme structure. Over it we have a family $\wh\cY\to \oZ$ of 
stable pairs $(\hY, \sum_{i=0,1,2} \frac12(R_i + G_i + B_i))$.

Stable pairs are described by partial matroid tilings of $\Delta\ubur$
covering the window $\Delta\halfwts(3,9)$, which can be computed
explicitly. Alternatively, one can start by classifying the tilings of
the small window $\Delta\halfwts(3,9)$ itself, a significantly easier
task. The latter tilings describe the irreducible components $\hY_k$ of
$\hY$. The possibilities for the incidence relations between the
curves $R_i$, $G_i$, $B_i$ can then be added in the second step.

The points $p_R$, $p_G$, $p_B$ are log centers of the pair
$(\bP^2, \sum\frac12(R_i+G_i+B_i)$: on the blowup the exceptional
divisors $R_3$, $G_3$, $B_3$ have discrepancy~$-1$. One checks that in
the family $\wh\cY$ these points give three disjoint sections and that
every fiber is smooth at these points. (This is a special case of a
general phenomenon.) Let $\wt\cY\to\wh\cY$ be the blowup at these
sections, and consider the divisor
\begin{displaymath}
  K_{\wt\cY} + \sum_{i=0}^3 \tfrac12 (\cR_i+\cG_i+\cB_i) = 
  f^*\big( K_{\wh\cY} + \sum_{i=0}^2 \tfrac12 (\cR_i+\cG_i+\cB_i) \big)
  - \tfrac12 (\cR_3+\cG_3+\cB_3).
\end{displaymath}
One checks that this divisor is big and nef on each fiber. By the
Basepoint-Free Theorem \cite[Thm.~3.24]{kollar1998birational-geometry}
it is relatively semiample and defines a contraction to a family
$\cY\to\oZ$ of stable pairs $(Y,\sum_{i=0}^3
\frac12(R_i+G_i+B_i)$. Some fibers in this family may be isomorphic. We
have a classifying map to the moduli space of stable pairs. Let
\begin{displaymath} 
  \oZ \to \oZ' \to \oM\slc
\end{displaymath}
be its Stein factorization. Then $\oZ'\to\oM\slc$ is a finite
birational map to the closure of $Z$ in $\oM\slc$. The normalization
of $\oZ'$ provides the required compactification for the moduli of the
pairs $\big(\Sigma, \sum(\frac12 R_i+\frac12 G_i+\frac12 B_i)\big)$.

\begin{remark}
  It follows that $\oM(\frac12)$ is the normalization of $\oZ'$.
\end{remark}

\begin{remark}
  If one is interested only in the irreducible components $Y_k$ of $Y$
  then some additional considerations show that it suffices to look only
  at the tilings of the polytope $\Delta\ubur\halfwts$ itself, instead
  of computing partial covers of $\Delta(3,9)$ containing
  $\Delta\ubur\halfwts$. In other words, the grayed out inequalities
  can be ignored.  This makes the computation much easier.
\end{remark}

\begin{remark}
  The same computations work for the window 
  $$0\le r_0,r_3,g_0,g_3,b_0,b_3 \le\tfrac12,\quad
  0\le r_1,r_2,g_1,g_2,b_1,b_2\le c$$ for any $\frac14<c\le
  \frac12$. For $\frac14<c\le\frac13$, $\frac13<c\le\frac25$ and
  $\frac25<c<\frac12$ respectively they produce the spaces
  $\oM(\frac13)$, $\oM(\frac25)$, $\oM(\frac12^-)$ of
  Section~\ref{sec:25}. For smaller $c$ the polytopes 6 and 8 of
  Table~\ref{tab:Y-irred-comps} do not intersect $\Delta\ubur\halfwts$,
  so some tilings become simpler. By the general theory
    \cite{alexeev2015moduli-weighted} there are reduction
  morphisms between the moduli spaces
  $\oM(\frac13)\gets \oM(\frac25) \gets \oM(\frac12^-)=\oM(\frac12)$,
  same as in Section~\ref{sec:25}. 
\end{remark}

\ifshort
\else
\fi

\bibliographystyle{amsalpha} 

\begin{thebibliography}{MLPR09}

\bibitem[AB21]{ascher2021moduli-of-weighted}
Kenneth Ascher and Dori Bejleri, \emph{Moduli of weighted stable elliptic
  surfaces and invariance of log plurigenera}, Proc. Lond. Math. Soc. (3)
  \textbf{122} (2021), no.~5, 617--677, With an appendix by Giovanni
  Inchiostro. \MR{4258169}

\bibitem[ABE22]{alexeev2022compactifications-moduli}
Valery Alexeev, Adrian Brunyate, and Philip Engel, \emph{Compactifications of
  moduli of elliptic {K}3 surfaces: {S}table pair and toroidal}, Geom. Topol.
  \textbf{26} (2022), no.~8, 3525--3588. \MR{4562567}

\bibitem[AE22]{alexeev2022mirror-symmetric}
Valery Alexeev and Philip Engel, \emph{Compactifications of moduli spaces of
  {K3} surfaces with a nonsymplectic involution}, arXiv:2208.10383 (2022).

\bibitem[AE23]{alexeev2023compact}
\bysame, \emph{Compact moduli of {K}3 surfaces}, Ann. of Math. (2) \textbf{198}
  (2023), no.~2, 727--789. \MR{4635303}

\bibitem[AET23]{alexeev2023stable-pair}
Valery Alexeev, Philip Engel, and Alan Thompson, \emph{Stable pair
  compactification of moduli of {K}3 surfaces of degree 2}, J. Reine Angew.
  Math. \textbf{799} (2023), 1--56. \MR{4595306}

\bibitem[AH23]{alexeev2023secondary-burniat}
Valery Alexeev and Xiaoyan Hu, \emph{Explicit compactifications of moduli
  spaces of secondary {B}urniat surfaces}, Preprint (2023), arXiv:2309.11397.

\bibitem[AK10]{alexeev2010complete-moduli}
Valery Alexeev and Allen Knutson, \emph{{C}omplete moduli spaces of
  branchvarieties}, J. Reine Angew. Math. \textbf{639} (2010), 39--71, arXiv:
  math.AG/0602626. \MR{2608190}

\bibitem[Ale94]{alexeev1994boundedness-and-ksp-2}
Valery Alexeev, \emph{Boundedness and {$K\sp 2$} for log surfaces}, Internat.
  J. Math. \textbf{5} (1994), no.~6, 779--810. \MR{MR1298994 (95k:14048)}

\bibitem[Ale02]{alexeev2002complete-moduli}
\bysame, \emph{{C}omplete moduli in the presence of semiabelian group action},
  Ann. of Math. (2) \textbf{155} (2002), no.~3, 611--708. \MR{2003g:14059}

\bibitem[Ale15]{alexeev2015moduli-weighted}
\bysame, \emph{Moduli of weighted hyperplane arrangements}, Advanced Courses in
  Mathematics. CRM Barcelona, Birkh\"{a}user/Springer, Basel, 2015, Edited by
  Gilberto Bini, Mart\'{i} Lahoz, Emanuele Macr\`\i and Paolo Stellari.
  \MR{3380944}

\bibitem[Ale23]{alexeev2023kappa-classes}
\bysame, \emph{Kappa classes on {KSBA} spaces}, Preprint (2023),
  arXiv:2309.14842.

\bibitem[AP12]{alexeev2012non-normal-abelian}
Valery Alexeev and Rita Pardini, \emph{Non-normal abelian covers}, Compos.
  Math. \textbf{148} (2012), no.~4, 1051--1084. \MR{2956036}

\bibitem[AP23]{alexeev2009explicit-compactifications}
\bysame, \emph{Explicit compactifications of moduli spaces of {C}ampedelli and
  {B}urniat surfaces}, Preprint (2023), arXiv:0901.4431v3.

\bibitem[BS92]{billera1992fiber-polytopes}
Louis~J. Billera and Bernd Sturmfels, \emph{Fiber polytopes}, Ann. of Math. (2)
  \textbf{135} (1992), no.~3, 527--549. \MR{1166643 (93e:52019)}

\bibitem[Cam32]{campedelli1932sopra-alcuni}
L.~Campedelli, \emph{Sopra alcuni piani doppi notevoli con curva di diramazioni
  del decimo ordine}, Atti Acad. Naz. Lincei \textbf{15} (1932), 536--542.

\bibitem[CFP{\etalchar{+}}23]{coughlan2023Tdivisors}
Stephen Coughlan, Marco Franciosi, Rita Pardini, Julie Rana, and S\"{o}nke
  Rollenske, \emph{On {T}-divisors and intersections in the moduli space of
  stable surfaces {$\overline{\mathfrak M}_{1,3}$}}, J. Lond. Math. Soc. (2)
  \textbf{107} (2023), no.~2, 750--776. \MR{4549144}

\bibitem[DO88]{dolgachev1988point-sets}
Igor Dolgachev and David Ortland, \emph{Point sets in projective spaces and
  theta functions}, Ast\'{e}risque (1988), no.~165, 210. \MR{1007155}

\bibitem[FFP22]{Fantechi2022semismooth}
Barbara Fantechi, Marco Franciosi, and Rita Pardini, \emph{Smoothing
  semi-smooth stable {G}odeaux surfaces}, Algebr. Geom. \textbf{9} (2022),
  no.~4, 502--512. \MR{4450623}

\bibitem[FPR17]{Franciosi2017K21}
Marco Franciosi, Rita Pardini, and S\"{o}nke Rollenske, \emph{Gorenstein stable
  surfaces with {$K^2_X=1$} and {$p_g>0$}}, Math. Nachr. \textbf{290} (2017),
  no.~5-6, 794--814. \MR{3636379}

\bibitem[FPRR22]{Franciosi2022Tsing}
Marco Franciosi, Rita Pardini, Julie Rana, and S\"{o}nke Rollenske,
  \emph{I-surfaces with one {T}-singularity}, Boll. Unione Mat. Ital.
  \textbf{15} (2022), no.~1-2, 173--190. \MR{4390548}

\bibitem[GPSZ22]{gallardo2022unimodal}
Patricio Gallardo, Gregory Pearlstein, Luca Schaffler, and Zheng Zhang,
  \emph{Unimodal singularities and boundary divisors in the {KSBA} moduli of a
  class of {H}orikawa surfaces}, 2022, arXiv:2209.08877.

\bibitem[Has03]{hassett2003moduli-spaces}
Brendan Hassett, \emph{{M}oduli spaces of weighted pointed stable curves}, Adv.
  Math. \textbf{173} (2003), no.~2, 316--352. \MR{MR1957831 (2004b:14040)}

\bibitem[HKT06]{hacking2006compactification-moduli}
Paul Hacking, Sean Keel, and Jenia Tevelev, \emph{Compactification of the
  moduli space of hyperplane arrangements}, J. Algebraic Geom. \textbf{15}
  (2006), no.~4, 657--680. \MR{2237265}

\bibitem[HS04]{haiman2004multigraded-hilbert}
Mark Haiman and Bernd Sturmfels, \emph{{M}ultigraded {H}ilbert schemes}, J.
  Algebraic Geom. \textbf{13} (2004), no.~4, 725--769. \MR{2073194
  (2005d:14006)}

\bibitem[Hu14]{hu14compactifications-of-moduli}
Xiaoyan Hu, \emph{The compactifications of moduli spaces of {B}urniat surfaces
  with {$2\leq K^{2} \leq5$}}, Ph.D. thesis, University of Georgia, 2014.

\bibitem[Inc20]{inchiostro20moduli-weierstrass}
Giovanni Inchiostro, \emph{Moduli of {W}eierstrass fibrations with marked
  section}, Adv. Math. \textbf{375} (2020), 107374, 57. \MR{4137071}

\bibitem[Ino94]{inoue1994some-new-surfaces}
Masahisa Inoue, \emph{Some new surfaces of general type}, Tokyo J. Math.
  \textbf{17} (1994), no.~2, 295--319. \MR{MR1305801 (95j:14048)}

\bibitem[Kat89]{kato1989logarithmic-structures}
Kazuya Kato, \emph{Logarithmic structures of {F}ontaine-{I}llusie}, Algebraic
  analysis, geometry, and number theory ({B}altimore, {MD}, 1988), Johns
  Hopkins Univ. Press, Baltimore, MD, 1989, pp.~191--224. \MR{1463703}
 
\bibitem[KK10]{kollar2010-dubois}
J\'{a}nos Koll\'{a}r and S\'{a}ndor~J. Kov\'{a}cs, \emph{Log canonical
  singularities are {D}u {B}ois}, J. Amer. Math. Soc. \textbf{23} (2010),
  no.~3, 791--813. \MR{2629988}


\bibitem[KM98]{kollar1998birational-geometry}
J\'{a}nos Koll\'{a}r and Shigefumi Mori, \emph{Birational geometry of algebraic
  varieties}, Cambridge Tracts in Mathematics, vol. 134, Cambridge University
  Press, Cambridge, 1998, With the collaboration of C. H. Clemens and A. Corti,
  Translated from the 1998 Japanese original. \MR{1658959}

 
  \bibitem[Kol13]{kollar2013singularities}
J\'{a}nos Koll\'{a}r, \emph{Singularities of the minimal model program. With the collaboration of S\'andor Kov\'acs},
Cambridge Tracts in Mathematics, vol.  200, Cambridge University Press Cambridge, 2023.

\bibitem[Kol23]{kollar2023families-of-varieties}
J\'{a}nos Koll\'{a}r, \emph{Families of varieties of general type}, Cambridge
  Tracts in Mathematics, vol. 231, Cambridge University Press, Cambridge, 2023.
  \MR{4566297}

\bibitem[KSB88]{kollar1988threefolds-and-deformations}
J.~Koll{\'a}r and N.~I. Shepherd-Barron, \emph{{T}hreefolds and deformations of
  surface singularities}, Invent. Math. \textbf{91} (1988), no.~2, 299--338.
  \MR{MR922803 (88m:14022)}

\bibitem[KSZ91]{kapranov1991quotients-toric}
M.~M. Kapranov, B.~Sturmfels, and A.~V. Zelevinsky, \emph{Quotients of toric
  varieties}, Math. Ann. \textbf{290} (1991), no.~4, 643--655. \MR{1119943}

\bibitem[Liu12]{liu2012stable-degenerations}
Wenfei Liu, \emph{Stable degenerations of surfaces isogenous to a product
  {II}}, Trans. Amer. Math. Soc. \textbf{364} (2012), no.~5, 2411--2427.
  \MR{2888212}

\bibitem[LR14]{liu-rollenske2014}
Wenfei Liu and S{\"o}nke Rollenske, \emph{Pluricanonical maps of stable log
  surfaces}, Adv. Math. \textbf{258} (2014), 69--126 (English).

\bibitem[MFK94]{mumford1994geometric-invariant}
D.~Mumford, J.~Fogarty, and F.~Kirwan, \emph{Geometric invariant theory}, third
  ed., Ergebnisse der Mathematik und ihrer Grenzgebiete (2) [Results in
  Mathematics and Related Areas (2)], vol.~34, Springer-Verlag, Berlin, 1994.
  \MR{MR1304906 (95m:14012)}

\bibitem[Miy77]{miyaoka1977on-numerical}
Y.~Miyaoka, \emph{On numerical {C}ampedelli surfaces}, Complex analysis and
  algebraic geometry, Iwanami Shoten, Tokyo, 1977, pp.~113--118. \MR{0447258}

\bibitem[MLP01]{lopes2001connected-component}
Margarida~Mendes~Lopes and Rita~Pardini, \emph{{A} connected component of the moduli
  space of surfaces with {$p\sb g=0$}}, Topology \textbf{40} (2001), no.~5,
  977--991.

\bibitem[MLPR09]{lopes2009campedelli-surfaces}
Margarida Mendes~Lopes, Rita Pardini, and Miles Reid, \emph{Campedelli surfaces
  with fundamental group of order 8}, Geom. Dedicata \textbf{139} (2009),
  49--55. \MR{2481836}

\bibitem[Mol21]{molcho2021universl-stacky}
Sam Molcho, \emph{Universal stacky semistable reduction}, Israel J. Math.
  \textbf{242} (2021), no.~1, 55--82. \MR{4282076}

\bibitem[Par91]{pardini1991abelian-covers}
Rita Pardini, \emph{{A}belian covers of algebraic varieties}, J. Reine Angew.
  Math. \textbf{417} (1991), 191--213. \MR{1103912 (92g:14012)}

\bibitem[PS02]{peeva2002toric-hilbert}
Irena Peeva and Mike Stillman, \emph{{T}oric {H}ilbert schemes}, Duke Math. J.
  \textbf{111} (2002), no.~3, 419--449. \MR{1885827 (2003m:14008)}

\bibitem[Rei83]{reid1983minimal-models}
Miles Reid, \emph{Minimal models of canonical {$3$}-folds}, Algebraic varieties
  and analytic varieties ({T}okyo, 1981), Adv. Stud. Pure Math., vol.~1,
  North-Holland, Amsterdam, 1983, pp.~131--180. \MR{715649}

\bibitem[Rol10]{rollenske2010compact-moduli}
S\"{o}nke Rollenske, \emph{Compact moduli for certain {K}odaira fibrations},
  Ann. Sc. Norm. Super. Pisa Cl. Sci. (5) \textbf{9} (2010), no.~4, 851--874.
  \MR{2789478}

\bibitem[{Sag}22]{sagemath}
{Sage Developers}, \emph{{S}agemath, the {S}age {M}athematics {S}oftware
  {S}ystem ({V}ersion 9.5)}, 2022, {\tt https://www.sagemath.org}.

\bibitem[vO06]{vanopstall2006stable-degenerations1}
Michael van Opstall, \emph{Stable degenerations of surfaces isogenous to a
  product of curves}, Proc. Amer. Math. Soc. \textbf{134} (2006), no.~10,
  2801--2806. \MR{2231601}

\end{thebibliography}

\newcommand{\etalchar}[1]{$^{#1}$}
\def\cprime{$'$}
\providecommand{\bysame}{\leavevmode\hbox to3em{\hrulefill}\thinspace}
\providecommand{\MR}{\relax\ifhmode\unskip\space\fi MR }
\providecommand{\MRhref}[2]{%
  \href{http://www.ams.org/mathscinet-getitem?mr=#1}{#2}
}
\providecommand{\href}[2]{#2}

\ifshort
\else

\fi

\end{document}